\documentclass{article}%
\usepackage{amsmath}%
\usepackage{amsfonts}%
\usepackage{amssymb}%
\usepackage{graphicx}
\newtheorem{theorem}{Theorem}

\newtheorem{corollary}[theorem]{Corollary}

\newtheorem{definition}[theorem]{Definition}

\newtheorem{lemma}[theorem]{Lemma}

\newtheorem{proposition}[theorem]{Proposition}

\newenvironment{proof}[1][Proof]{\noindent\textbf{#1.} }{\ \rule{0.5em}{0.5em}}
\begin{document}

\title{A bracket polynomial for graphs. III. \\Vertex weights}
\author{Lorenzo Traldi\\Lafayette College\\Easton, Pennsylvania 18042}
\date{}
\maketitle

\begin{abstract}
In earlier work the Kauffman bracket polynomial was extended to an invariant
of marked graphs, i.e., looped graphs whose vertices have been partitioned
into two classes (marked and not marked). The marked-graph bracket polynomial
is readily modified to handle graphs with weighted vertices. We present
formulas that simplify the computation of this weighted bracket for graphs
that contain twin vertices or are constructed using graph composition, and we
show that graph composition corresponds to the construction of a link diagram
from tangles.

\bigskip


\bigskip


\end{abstract}

\section{Introduction}

In this paper a \textit{graph} $G=(V(G),E(G))$ may have loops or multiple
edges; it may also contain free loops, which are connected components that
contain neither vertices nor edges.\ A \textit{marked graph} is a graph whose
vertex-set has been partitioned into two subsets, either of which may be
empty. The vertices in one cell of the partition are \textit{unmarked} and the
vertices in the other cell are \textit{marked}; in figures we indicate marked
vertices with the letter \textit{c}. If $G$ is a marked graph with
$V(G)=\{v_{1},...,v_{n}\}$ then the \textit{Boolean adjacency matrix}
$\mathcal{A}(G)$ is an $n\times n$ matrix over the two-element field $GF(2)$,
with entries $\mathcal{A}(G)_{ii}=1$ if $v_{i}$ is looped and for $i\neq j$,
$\mathcal{A}(G)_{ij}=1$ if $v_{i}$ and $v_{j}$ are adjacent. For $T\subseteq
V(G)$ we denote by $\mathcal{A}(G)_{T}$ the matrix obtained from
$\mathcal{A}(G)$ by first changing the $i^{th}$ diagonal entry whenever
$v_{i}\in T$, and then removing the $i^{th}$ row and column whenever $v_{i}$
is marked and the $i^{th}$ diagonal entry is 0. The \textit{marked-graph
bracket polynomial} of $G$ is defined by the formula
\[
\lbrack G]=d^{\phi}\cdot\sum_{T\subseteq V(G)}A^{n-\left\vert T\right\vert
}B^{\left\vert T\right\vert }d^{\nu(\mathcal{A}(G)_{T})},
\]
where $\phi$ is the number of free loops in $G$ and $\nu$ is the
$GF(2)$-nullity of a matrix \cite{T1, TZ}.%

\begin{figure}
[ptb]
\begin{center}
\includegraphics[
trim=0.974160in 8.751122in 0.941990in 0.938346in,
height=0.787in,
width=4.7738in
]%
{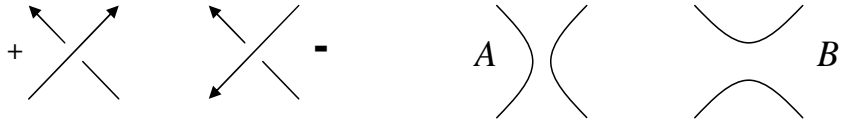}%
\caption{Signs and smoothings of classical crossings.}%
\label{flyfig1}%
\end{center}
\end{figure}

An \textit{oriented regular link diagram} $D$ consists of oriented, piecewise
smooth, closed curves in the plane; the curves intersect (and self-intersect)
only at a finite number of transverse \textit{crossings}. There are two kinds
of crossings, classical and virtual, with underpassing and overpassing arcs
specified at the classical crossings. A classical crossing has two
\textit{smoothings}, one denoted $A$ and the other denoted $B$, as in Figure
\ref{flyfig1}. If $D$ has $n$ classical crossings then it has $2^{n}$
\textit{Kauffman states}, obtained by applying either the $A$ or the $B$
smoothing at each classical crossing. Given a state $S$ let $a(S)$ denote the
number of $A$ smoothings in $S$, $b(S)=n-a(S)$ the number of $B$ smoothings in
$S$, and $c(S)$ the number of closed curves in $S$, including any
crossing-free components that might appear in $D$. Then the (three-variable)
\textit{Kauffman bracket polynomial} of $D$ \cite{Kau, Kv, Kd} is
\[
\lbrack D]=\sum_{S}A^{a(S)}B^{b(S)}d^{c(S)-1}.
\]
The fact that $A$, $B$ and $d$ are independent variables is indicated by using
$[D]$ rather than the more familiar notation $\left\langle D\right\rangle $.

\begin{figure}
[h]
\begin{center}
\includegraphics[
trim=1.106613in 3.344921in 0.536815in 1.079317in,
height=4.7323in,
width=4.9813in
]%
{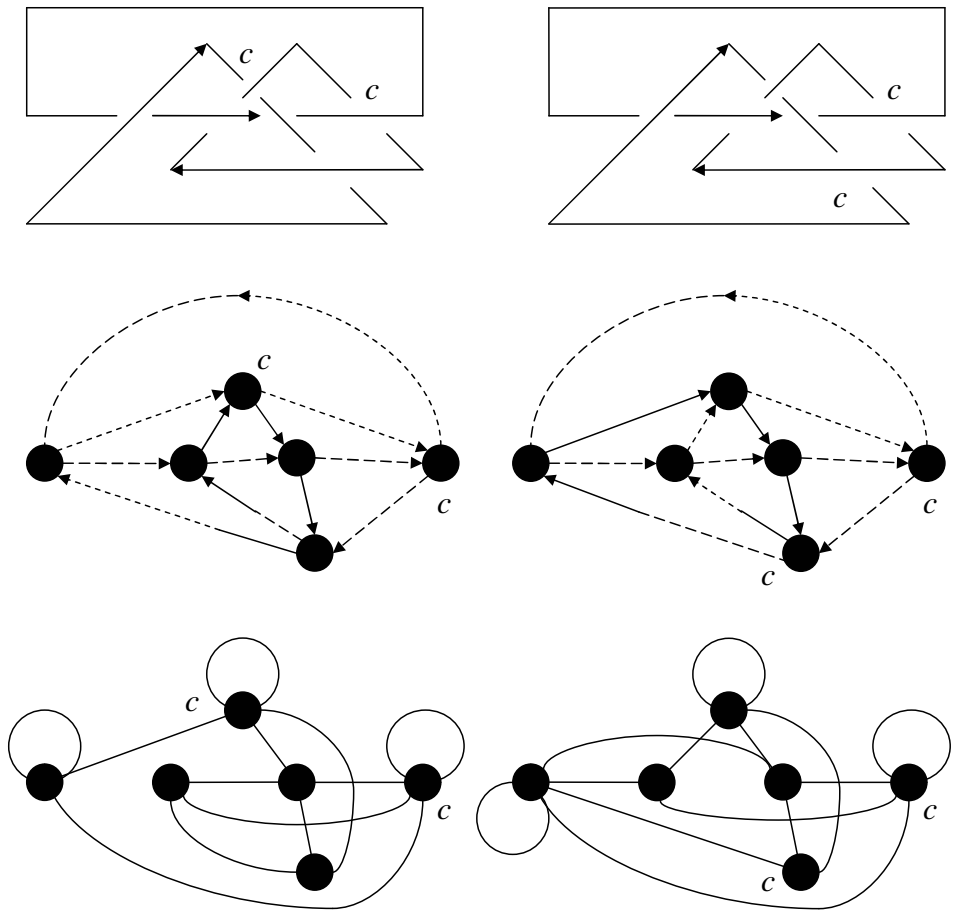}%
\caption{Two marked diagrams of the Borromean rings, the corresponding
directed universes with Euler circuits indicated by dashes, and the associated
looped interlacement graphs.}%
\label{borr3}%
\end{center}
\end{figure}

These two kinds of bracket polynomials are closely related. A link diagram $D
$ has an associated \textit{directed universe} $\vec{U}$, a 2-in, 2-out
digraph whose vertices correspond to the classical crossings of $D$ and whose
edges correspond to the arcs of $D$. ($\vec{U}$ also contains a free loop for
each link component that is crossing-free in $D$.) The undirected version of
$\vec{U}$ is denoted $U$. Let $C$ be a directed Euler system for $\vec{U}$,
i.e., a set containing one directed Euler circuit for each connected component
of $U$; $C$ must also contain every free loop of $U$. $C$ is completely
determined by specifying the classical crossings at which it does not follow
the incident link component(s). We say such crossings are \textit{marked}, and
we indicate them in figures with the letter \textit{c}. The \textit{looped
interlacement graph} $\mathcal{L}(D,C)$ \cite{T1, TZ} is a marked graph with a
vertex for each classical crossing of $D$; marked vertices correspond to
marked crossings, and looped vertices correspond to negative crossings. Two
vertices $v$ and $w$ are adjacent in $\mathcal{L}(D,C)$ if and only if they
are \textit{interlaced} with respect to $C$, i.e., there is a circuit of $C$
on which they appear in the order $v...w...v...w$ \cite{RR}. Different choices
of $C$ may give rise to different graphs $\mathcal{L}(D,C)$, as in\ Figure
\ref{borr3}, but every looped interlacement graph has marked-graph bracket
$[\mathcal{L}(D,C)]$ equal to the Kauffman bracket $[D]$.

The identity $[\mathcal{L}(D,C)]=[D]$ is derived from an equality connecting
nullities of matrices over $GF(2)$ with circuits in 4-regular graphs: if $S$
is a Kauffman state in a link diagram $D$ and $T\subseteq V(\mathcal{L}(D,C))
$ consists of those vertices where $S$ involves the $B$ smoothing, then the
equality tells us that $c(S)-1=\phi+\nu(\mathcal{A}(\mathcal{L}(D,C))_{T})$.
Several similar equalities have been discovered independently over the years.
This form of the equality originated as a result about permutations due to
Cohn and Lempel \cite{CL}, but related results appeared much earlier, in
Brahana's study of curves on surfaces \cite{Br}. Later authors worked in
combinatorics \cite{Be, Ma} or classical knot theory \cite{Me, S, Z}. See
\cite{Tb} for a thorough account.%

For alternating classical link diagrams, the equality connecting nullity and
circuits appears implicitly in the relationship between the Kauffman bracket
and certain combinatorial invariants: the Tutte polynomial \cite{Th}, Jaeger's
transition polynomial \cite{E, J}, or the interlace polynomial \cite{A2}. The
Tutte polynomial involves the checkerboard graph, and the latter two
polynomials involve a directed graph obtained from the universe of a link
diagram by reversing the orientation of every second edge. These relationships
also hold for non-alternating classical diagrams, if weighted edges or
vertices are used to distinguish between positive and negative classical
crossings \cite{K1, M, SW, T2, Tw}. Since Kauffman introduced virtual knot
theory \cite{Kv}, the combinatorial theory associated with the Kauffman
bracket of classical link diagrams has been extended to virtual link diagrams
in several ways. For checkerboard-colorable virtual links, little modification
of the classical theory is needed \cite{Ka}. For general virtual links there
is not a direct analogue of the checkerboard graph or the edge-reversed
universe graph, but the standard Tutte polynomial may be replaced by a
``topological''\ Tutte polynomial \cite{BR1, BR2, CP, CV}, or by a
``relative''\ or ``ported''\ Tutte polynomial \cite{Ch, DH}. Kauffman's
bracket may also be analyzed combinatorially without using any version of the
Tutte polynomial; Ilyutko and Manturov have developed a geometric approach
involving atoms and rotating circuits \cite{IM}, and Zulli and the present
author have developed a graph-theoretic approach using Euler circuits,
interlacement and marked graphs \cite{T1, TZ}.

Marked graphs are subject to two important equivalence relations, which
generalize natural equivalences between looped interlacement graphs. The finer
of the two equivalence relations generalizes the equivalence between looped
interlacement graphs $\mathcal{L}(D,C_{1})$ and $\mathcal{L}(D,C_{2})$
obtained from a single link diagram $D$. This equivalence relation is
generated by a graph-theoretic operation we call a \textit{marked pivot}; it
is a modified version of the pivot operation that describes the effect on
interlacement graphs of changing Euler systems in 2-in, 2-out digraphs
\cite{A2, K}. Marked pivots preserve the 3-variable marked-graph bracket
polynomial. The coarser of the two equivalence relations generalizes the
equivalence between looped interlacement graphs $\mathcal{L}(D_{1},C_{1})$ and
$\mathcal{L}(D_{2},C_{2})$ obtained from different diagrams of the same link
type. This equivalence relation is generated by graph-theoretic versions of
the Reidemeister moves. Marked-graph Reidemeister moves change the 3-variable
marked-graph bracket, but they preserve the marked-graph analogue of the Jones
polynomial, which is obtained (as usual) by replacing $A$ with $t^{-1/4}$,
replacing $B$ with $t^{1/4}$, replacing $d$ with $-t^{1/2}-t^{-1/2}$, and
multiplying by a suitable factor.

These two equivalence relations are certainly important, for they underlie the
knot-theoretic significance of the marked-graph bracket. Nevertheless we pay
little attention to them in this paper, because our purpose is not to discuss
the relationship between bracket polynomials of different graphs, but rather
to discuss the efficient computation of the bracket polynomial of a given graph.

The formulas that define the marked-graph bracket and the three-variable
Kauffman bracket are quite similar, so it may be surprising that recursive
descriptions of the two bracket polynomials are quite different. The Kauffman
bracket of a link diagram $D$ can be calculated by repeatedly applying a
single recursive step: choose a classical crossing in $D$, let $D_{A}$ and
$D_{B}$ be the two diagrams obtained by smoothing that crossing, and use the
formula $[D]=A[D_{A}]+B[D_{B}]$. The recursive description of the marked-graph
bracket given in \cite{T1, TZ} is considerably more complicated; four
different recursive steps are used, in different situations. For instance, one
step removes a loop on an unmarked vertex, at the cost of replacing the graph
in question with two smaller graphs; this particular step is related to the
Jones polynomial's fundamental identity $tV_{L^{-}}=t^{-1}V_{L^{+}}%
-(t^{1/2}-t^{-1/2})V_{L}$ \cite{Jo}.

If we rewrite the definition of the marked-graph bracket polynomial as
\[
\lbrack G]=d^{\phi}\cdot\sum_{T\subseteq V(G)}(\prod_{v\notin T}A)(\prod_{t\in
T}B)d^{\nu(\mathcal{A}(G)_{T})},
\]
then the variables $A$ and $B$ appear as vertex weights, similar to those that
appear in many combinatorial contexts ranging from electrical circuit theory
to statistical mechanics. (For instance, if $A=B=\frac{1}{2}$ then $[G]$ gives
the expected value of $d^{\phi+\nu(\mathcal{A}(G)_{T})}$ under the presumption
that $T$ is chosen by tossing a fair coin $n$ times, the $i^{th}$ toss
deciding whether or not $v_{i}\in T$.) The following generalization suggests itself.

\begin{definition}
\label{bracweigh} Suppose $G$ is a weighted, marked graph, i.e., a marked
graph given with functions $\alpha$ and $\beta$ mapping $V(G)$ into some
commutative ring $R$. Then the \emph{weighted marked-graph bracket polynomial}
of $G$ is
\[
\lbrack G]=d^{\phi}\cdot\sum_{T\subseteq V(G)}(\prod_{v\notin T}%
\alpha(v))(\prod_{t\in T}\beta(t))d^{\nu(\mathcal{A}(G)_{T})}%
\]

\end{definition}

If $v\in V(G)$ has $\alpha(v)=A$ and $\beta(v)=B$ then we say $v$ has
\emph{standard weights}. If $D$ is a link diagram then the vertices of
$G=\mathcal{L}(D,C)$ correspond to the classical crossings of $D$, so we may
think of $\alpha$ and $\beta$ as giving weights for the classical crossings of
$D$.

The recursive description of the marked-graph bracket polynomial given in
\cite{T1} extends directly to the weighted version of the polynomial. Vertex
weights may be used to make the recursion more efficient in several ways. The
most obvious simplification involves the recursive step mentioned above, used
to eliminate loops on unmarked vertices; it is completely unnecessary.

\begin{theorem}
\label{loops} Suppose $G$ is a weighed, marked graph. Let $G^{\prime}$ be the
graph obtained by removing every loop from $G$, and reversing the $\alpha$ and
$\beta$ weights of every looped vertex of $G$. Then $[G]=[G^{\prime}]$.
\end{theorem}

The value of Theorem \ref{loops} is easy to see: each time we use the theorem
instead of a loop-removing recursive step, the recursion proceeds with only
one graph to process rather than two.

The weighted marked-graph bracket polynomial also satisfies several analogues
of the series-parallel reductions of electrical circuit theory. These are
operations which consolidate certain vertices without changing the value of
the bracket. Here is one of them.

\begin{theorem}
\label{parallels} Suppose $v_{1},...,v_{k}$ are unlooped twins that form a
clique in $G$, i.e., $v_{1},...,v_{k}$ all have the same neighbors outside
$\{v_{1},...,v_{k}\}$ and they are all adjacent to each other. Let
$\rho=\left\vert \{i|v_{i}\text{ is not marked}\}\right\vert .$ Let
$(G-v_{2}-...-v_{k})^{\prime}$ be the graph obtained from $G-v_{2}-...-v_{k}$
by (i) marking $v_{1}$ if and only if $\rho$ is even, and (ii) changing the
weights of $v_{1}$ to
\[
\alpha^{\prime}(v_{1})=\prod_{i=1}^{k}\alpha(v_{i})\text{ \ and \ }%
\beta^{\prime}(v_{1})=d^{-1}\left(  -\alpha^{\prime}(v_{1})+\prod_{i=1}%
^{k}(\alpha(v_{i})+d\beta(v_{i}))\right)  .
\]
Then $[G]=[(G-v_{2}-...-v_{k})^{\prime}]$.
\end{theorem}

The value of Theorem \ref{parallels} is not quite as obvious as that of
Theorem \ref{loops}. A first impression might be that we are simply replacing
$k$ vertices with one vertex, but this impression is imprecise because the
complicated values of $\alpha^{\prime}(v_{1})$ and $\beta^{\prime}(v_{1})$
given in the theorem may be inconvenient. The computational cost of this
inconvenience depends on the implementation of arithmetic operations in the
particular ring being used for $R$. For instance, a\ natural example is a ring
$R$ of polynomials in variables $\alpha_{1},...,\alpha_{n},\beta_{1}%
,...,\beta_{n}$, with the variables used as vertex-weights; arithmetic in this
ring is very expensive because each polynomial involves coefficients of many
different monomials. Nevertheless, Theorem \ref{parallels} clearly has the
potential to be of significant value in general.%

\begin{figure}
[h]
\begin{center}
\includegraphics[
trim=0.935921in 8.826013in 1.204741in 0.810825in,
height=0.8242in,
width=4.6095in
]%
{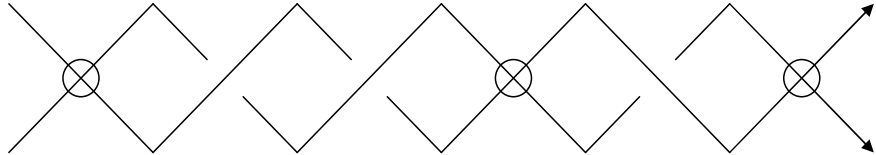}%
\caption{Twisted strands give rise to a clique of twin vertices in
$\mathcal{L}(D,C)$.}%
\label{twins2}%
\end{center}
\end{figure}

The most familiar situation in which Theorem \ref{parallels} arises involves
two coherently oriented strands of a link diagram, twisted around each other
to produce $k$ classical crossings and some number of virtual crossings.
Although the theorem specifies that the vertices are unlooped, negative
crossings may be handled simply by reversing their $\alpha$ and $\beta$
weights. An example with $k=3$ appears in Figure \ref{twins2}. (As usual, the
circled crossings are virtual.) If these twisted strands appear in a link
diagram $D$ with standard weights then the two negative crossings give rise to
unlooped vertices of $\mathcal{L}(D,C)$ with $\alpha=B$ and $\beta=A$, while
the positive crossing gives rise to an unlooped vertex with $\alpha=A$ and
$\beta=B$. Theorem \ref{parallels} tells us that the bracket is unchanged if
the three vertices of $\mathcal{L}(D,C)$ representing the portion of
$D$\ appearing in\ Figure \ref{twins2} are replaced with one unlooped,
unmarked vertex whose weights are $\alpha=AB^{2}$ and $\beta=(-AB^{2}%
+(A+Bd)(B+Ad)^{2})/d$. This new graph is the looped interlacement graph of a
diagram $D^{\prime}$ obtained by replacing the pictured portion of $D$ with a
single positive crossing carrying the indicated weights, and also a single
virtual crossing; the latter is needed whenever the total number of unmarked
classical and virtual crossings is even, to ensure that $C $ gives rise to an
Euler system of $D^{\prime}$.

In the classical case -- or more generally, the checkerboard-colorable case
\cite{Ka} -- twisting two strands around each other produces classical
crossings that give rise to series-parallel edges in the checkerboard graphs.
Theorem \ref{parallels} applies to all virtual link diagrams, not just the
checkerboard-colorable ones, but the extra generality comes at a price: the
theorem must be adjusted when there are marks on the vertices or
non-adjacencies among them, and not all of the adjusted versions are quite so
simple. For instance, if two oppositely oriented strands of a link diagram are
twisted around each other to produce a set of $k$ unmarked, nonadjacent twins
then a \textquotedblleft dual\textquotedblright\ of Theorem \ref{parallels}
requires that $k$ be odd. See Corollary \ref{parallelsdual}.

A third use for vertex weights is that graphs constructed from smaller graphs
using an appropriate version of Cunningham's \textit{composition} operation
\cite{Cu} have bracket polynomials that can be described by modifications of
the weights.%

\begin{figure}
[ptb]
\begin{center}
\includegraphics[
trim=1.130527in 7.886825in 0.672874in 0.804407in,
height=1.5316in,
width=4.8585in
]%
{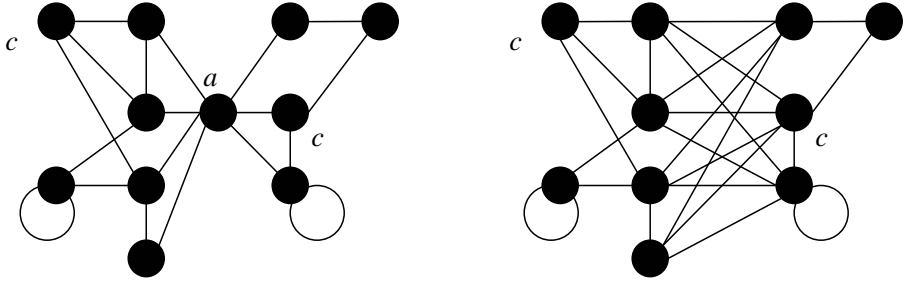}%
\caption{Composition of graphs.}%
\label{compfig}%
\end{center}
\end{figure}

\begin{definition}
\label{comp}A marked, weighted graph $G$ is the \emph{composition} of marked,
weighted graphs $F$ and $H$, $G=F\ast H$, if the following conditions hold.

(a) $V(F)\cap V(H)$ consists of a single unlooped, unmarked vertex $a$ that
has standard weights in both $F$ and $H$.

(b) The elements of $V(G)=V(F)\cup V(H)-\{a\}$ inherit their loops, marks and
weights from $F$ and $H$.

(c) $E(G)=E(F-a)\cup E(H-a)\cup\{\{v,w\}|\{v,a\}\in E(F)$ and $\{a,w\}\in
E(H)\}$.

(d) $F$ and $H$ do not share any free loop, and the free loops of $G$ are
those of $F$ and $H$.
\end{definition}

Requiring that $a$ have standard weights and be unlooped and unmarked ensures
that no significant information is lost when we remove $a$ in constructing
$F\ast H$. Note that Definition \ref{comp} includes the situation of Theorem
\ref{parallels}: if $F$ is a complete graph then $F-a$ is a clique of twins in
$F\ast H$.

The construction given in Definition \ref{comp}\ may seem to be merely a
technical notion from graph theory, but in Section 2 we show that it is
related to an important knot-theoretic idea: if a link diagram contains a
tangle, then the looped interlacement graph is a composition. Recall that a
subgraph of a graph $G$ is \textit{full} or \textit{induced} if it contains
every edge of $G$ incident on its vertices.

\begin{theorem}
\label{tangle} Suppose a link diagram $D$ contains a tangle, and let $C$ be
any directed Euler system for the universe of $D$. Then there are graphs $F$
and $H$ such that $\mathcal{L}(D,C)=F\ast H$, $F-a$ is the subgraph of
$\mathcal{L}(D,C)$ induced by the vertices corresponding to crossings outside
the tangle, and $H-a$ is the subgraph of $\mathcal{L}(D,C)$ induced by the
vertices corresponding to crossings inside the tangle. Moreover, it is
possible to choose $C$ so that at least one of $F,H$ has no marked vertex
adjacent to $a$.
\end{theorem}%

\begin{figure}
[tb]
\begin{center}
\includegraphics[
trim=0.660505in 4.546183in 1.071155in 1.075038in,
height=3.8346in,
width=4.913in
]%
{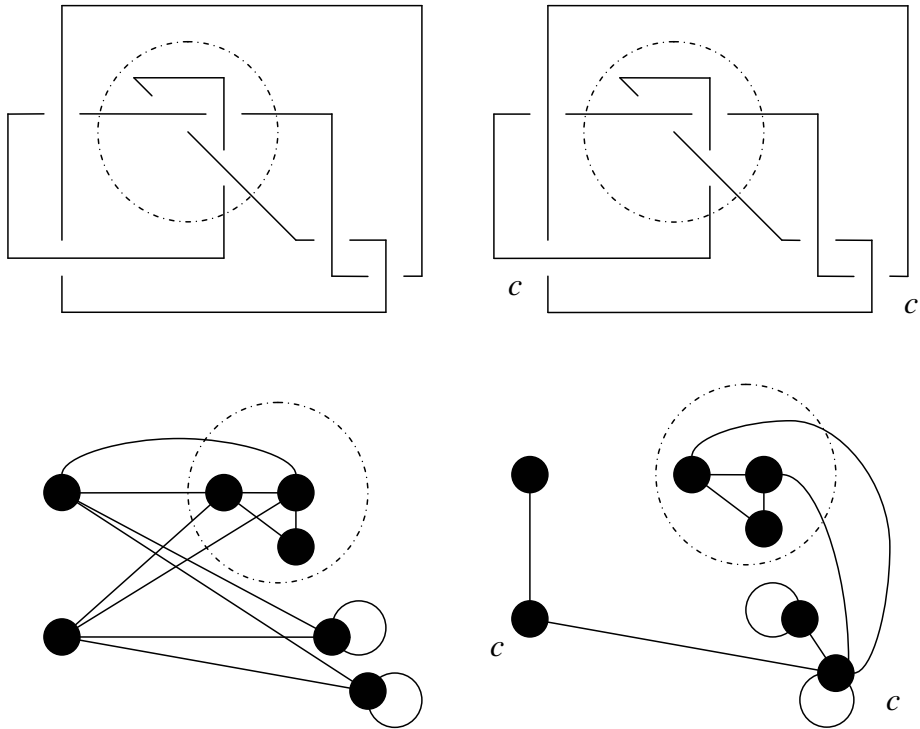}%
\caption{If $D$ contains a tangle then $\mathcal{L}(D,C)$ is a composition of
graphs.}%
\label{tanglex}%
\end{center}
\end{figure}

Two instances of Theorem \ref{tangle} are pictured in\ Figure \ref{tanglex}.
Observe that $D$ contains several tangles in addition to the one indicated by
the dashed circle. The interlacement graph on the left satisfies the last
sentence of the theorem for every tangle, as no vertex is marked. The
interlacement graph on the right, instead, satisfies the last sentence for
some tangles but not for others. For instance, it satisfies the last sentence
of the theorem for the indicated tangle, but not for the tangle that contains
the two crossings in the lower right-hand corner of the diagram.

For a fixed graph $F$, every composition $F\ast H$ is constructed from $H$ in
much the same way. In Section 5 we prove that similarly, every weighted
bracket polynomial $[F\ast H]$ is constructed from bracket polynomials
associated with $H$ in much the same way.

\begin{theorem}
\label{pjoin} Let $F$ be a marked, weighted graph with an unlooped, unmarked
vertex $a$ that has standard weights. Then there are weights $\alpha^{\prime
}(a)$, $\beta^{\prime}(a)$, and $\alpha^{\prime}(a_{m})$ that depend only on
$F$ and $a$, and have the following \textquotedblleft
universal\textquotedblright\ property: every composition $F\ast H$ obtained by
applying Definition \ref{comp} to a graph $H$ in which no neighbor of $a$ is
marked has
\[
\lbrack F\ast H]=[H^{\prime}]+[H_{m}^{\prime}],
\]
where $H^{\prime}$ is obtained from $H$ by changing the weights of $a$ to
$\alpha^{\prime}(a)$ and $\beta^{\prime}(a)$, and $H_{m}^{\prime}$ is obtained
from $H$ by marking $a$ and changing its weights to $\alpha^{\prime}(a_{m})$
and $\beta^{\prime}(a_{m})=0$.
\end{theorem}

Choosing to have $\beta^{\prime}(a_{m})=0$ in Theorem \ref{pjoin} is a matter
of convenience rather than necessity. The proof actually shows that the value
of $\beta^{\prime}(a_{m})$ is arbitrary, but the sum $\beta^{\prime}%
(a)+\beta^{\prime}(a_{m})$ must be correct. That is, for any $r\in R$ the
theorem still holds if we change $\beta^{\prime}(a_{m})$ from $0$ to $r$, and
also change $\beta^{\prime}(a)$ to $\beta^{\prime}(a)-r$. In particular,
Theorem \ref{pjoin} remains valid if $\beta^{\prime}(a)$ and $\beta^{\prime
}(a_{m})$ are interchanged.

It is hard to assess the loss of generality associated with Theorem
\ref{pjoin}'s hypothesis that $a$ have no marked neighbor in $H$, because
vertex marks do not have a special significance in general. For tangles in
link diagrams, though, the last sentence of Theorem \ref{tangle} tells us that
this hypothesis does not entail any loss of generality as long as we are
willing to reverse the roles of $F$ and $H$, i.e., to reverse the
\textquotedblleft inside\textquotedblright\ and \textquotedblleft
outside\textquotedblright\ of the tangle.

The idea that tangles are building blocks for classical link diagrams appears
in Conway's seminal paper \cite{C}. This idea leads naturally to the
observation that a recursively defined classical link invariant can be
calculated in a \textquotedblleft tangle-based\textquotedblright\ manner:
first eliminate all the crossings inside a particular tangle and then collect
like terms, before proceeding to eliminate the crossings inside another
tangle. See \cite{R} for a detailed analysis of the computational complexity
of such a tangle-based calculation of classical link invariants. The process
of building up link diagrams from tangles corresponds to the process of
building up checkerboard graphs as 2-sums of smaller graphs, and this
correspondence is useful for the Kauffman bracket because the Tutte polynomial
of a 2-sum can be described using the Tutte polynomials of the smaller graphs
\cite{Bry}. The correspondence between tangles and 2-sums was mentioned in
\cite{SW}, and its computational significance has recently been analyzed in
\cite{DEZ}. The result of a tangle-based calculation of the Kauffman bracket
of a classical diagram containing a tangle $T$ may be represented
schematically as $[T]=\gamma_{1}[\asymp]+\gamma_{2}[)(]$, where $\gamma_{1}$
and $\gamma_{2}$ are coefficients that result from the collection of terms.
Virtual crossings necessitate a third term; schematically, $[T]=\gamma
_{1}[\asymp]+\gamma_{2}[)(]+\gamma_{3}[\otimes]$. At first glance this
three-term schematic formula may seem different from the two-term formula of
Theorem \ref{pjoin}, but they are actually quite similar; each has three
degrees of freedom, represented by $\gamma_{1},\gamma_{2},\gamma_{3}$ in the
schematic formula and $\alpha(a),\beta(a),\alpha(a_{m})$ in Theorem
\ref{pjoin}.

Formulas for $\alpha(a)$, $\beta(a)$, and $\alpha(a_{m})$ are presented in
Corollaries \ref{pjoin1}, \ref{pjoin3} and \ref{pjoin4}. The first two are
obtained by using several specific graphs for $H$, and then solving the
resulting equations. The third breaks Definition \ref{bracweigh} into three
\textquotedblleft sub-sums.\textquotedblright

Suppose $F$ is a marked, weighted graph with an unlooped, unmarked vertex $a$
that has standard weights. For $i\neq j\in\{0,1\}$ let $F^{ij}$ be the graph
obtained from $F$ by replacing $a$ with a vertex $v_{ij}$ whose weights are
$\alpha(v_{ij})=i$ and $\beta(v_{ij})=j$. Then $F-a$, $F^{10}$ and $F^{01}$
are compositions $F\ast H^{0}$, $F\ast H^{10}$ and $F\ast H^{01} $
respectively, where $H^{0}$ is just $a$ and $H^{ij}$ has the two adjacent,
unlooped, unmarked vertices $a$ and $v_{ij}$. Theorem \ref{pjoin} gives three
equations.
\[
\lbrack F-a]=[F\ast H^{0}]=[H^{0\prime}]+[H_{m}^{0\prime}]=d\alpha
(a)+\beta(a)+\alpha(a_{m})
\]
\[
\lbrack F^{10}]=[F\ast H^{10}]=[H^{10\prime}]+[H_{m}^{10\prime}]=\alpha
(a)+\beta(a)+d\alpha(a_{m})
\]
\[
\lbrack F^{01}]=[F\ast H^{01}]=[H^{01\prime}]+[H_{m}^{01\prime}]=\alpha
(a)+d\beta(a)+\alpha(a_{m})
\]
We deduce the following.

\begin{corollary}
\label{pjoin1} The weights mentioned in Theorem \ref{pjoin} are given by these
formulas.
\begin{align*}
(2-d-d^{2})\alpha(a)  & =-(d+1)[F-a]+[F^{10}]+[F^{01}]\\
(2-d-d^{2})\beta(a)  & =[F-a]+[F^{10}]-(d+1)[F^{01}]\\
(2-d-d^{2})\alpha(a_{m})  & =[F-a]-(d+1)[F^{10}]+[F^{01}]
\end{align*}

\end{corollary}

To assess the computational significance of Theorem \ref{pjoin} and Corollary
\ref{pjoin1}, consider that the number of steps in an implementation of the
recursive algorithm for calculating $[G]$ is roughly $2^{\left\vert
V(G)\right\vert }$. (This is only a rough count rather than a precise
determination of computational complexity, because it ignores both the
computational cost of arithmetic in $R$ and the fact that different branches
of a calculation may require different numbers of steps.) Consequently the
number of steps involved in a direct computation of $[F\ast H]$ is roughly the
product of the numbers of steps involved in separate computations of $[F-a]$
and $[H-a]$. Corollary \ref{pjoin1} shows that the weights of $a$ and $a_{m}$
may be calculated using the bracket polynomials of $F-a$ and two graphs that
have $\left\vert V(F)\right\vert $ vertices apiece, so we may estimate the
number of steps involved in finding these weights as roughly five times the
number of steps involved in finding $[F-a]$. Theorem \ref{pjoin} then
expresses $[F\ast H]$ as the sum of the bracket polynomials of two $\left\vert
V(H)\right\vert $-vertex graphs; computing each of these brackets takes
roughly twice as many steps as computing $[H-a]$. All in all, Theorem
\ref{pjoin} and Corollary \ref{pjoin1} tell us that a rough upper bound on the
number of steps required to compute $[F\ast H]$ is on the order of five times
the sum of the numbers of steps required to calculate $[F-a]$ and $[H-a]$
separately. Five times the sum of two positive integers is generally
considerably smaller than the product of the two integers, so\ despite the
rough counting it is clear that when they are applicable, Theorem \ref{pjoin}
and Corollary \ref{pjoin1} can be considerably\ more efficient than direct computation.

Theorem \ref{pjoin} focuses on $F$, but compositions are symmetric and
consequently the theorem may be applied to both $F$ and $H$, so long as
neither contains a marked neighbor of $a$. Unlike the hypothesis of Theorem
\ref{pjoin} that only $H$ have no marked vertex adjacent to $a$, this double
hypothesis is a significant restriction even for link diagrams. The result is
still useful, though; for instance, every classical or virtual knot diagram
has an Euler system with respect to which there are no marked vertices at all.

\begin{corollary}
\label{pjoin2} Let $a(F)$ and $a_{m}(F)$ be the new weighted vertices
associated to $F$ and $a$ in Theorem \ref{pjoin}, and let $a(H)$ and
$a_{m}(H)$ be the new vertices associated in the same way to $H$ and $a$.
(That is, they are obtained by interchanging $F$ and $H$ in Theorem
\ref{pjoin}.) If neither $F$ nor $H$ contains a marked neighbor of $a$ then
\begin{gather*}
\lbrack F\ast H]=\alpha(a(F))\cdot(\alpha(a(H))+\beta(a(H)))+\beta
(a(F))\cdot(\alpha(a(H))+\beta(a(H))d)\\
+(\alpha(a(F))d+\beta(a(F)))\cdot\alpha(a_{m}(H))\\
+\alpha(a_{m}(F))\cdot(\alpha(a(H))d+\beta(a(H)))+\alpha(a_{m}(F))\cdot
\alpha(a_{m}(H)).
\end{gather*}

\end{corollary}

The corollary is proven as follows. Let $P_{1}$ be the graph with two adjacent
vertices $a$ and $b$, with $b$ an unlooped, unmarked vertex carrying the
weights of $a(F)$. Let $P_{1m}$ have two adjacent vertices $a$ and $b_{m}$
with $b_{m}$ an unlooped, marked vertex carrying the weights of $a_{m}(F)$.
Three applications of Theorem \ref{pjoin} tell us that
\begin{align*}
\lbrack F\ast H]  & =[H^{\prime}]+[H_{m}^{\prime}]=[H\ast P_{1}]+[H\ast
P_{1m}]\\
& =[P_{1}^{\prime}]+[(P_{1}^{\prime})_{m}]+[(P_{1m})^{\prime}]+[(P_{1m}%
)_{m}^{\prime}].
\end{align*}

Before proceeding we should express our gratitude to V. O. Manturov for his
comments and corrections regarding earlier versions of the paper.

\section{Diagrams and tangles}

Suppose $D$ is a diagram of a $\mu$-component link $L=K_{1}\cup...\cup K_{\mu
}$, and let $\Gamma$ denote the graph with $V(\Gamma)=\{v_{1},...,v_{\mu}\}$
and $E(\Gamma)=\{\{v_{i},v_{j}\}$ such that $D$ contains a classical crossing
involving $K_{i}$ and $K_{j}\}$. Suppose $C$ is a directed Euler system of
$\vec{U}$. As we traverse one of the circuits of $C$, we pass from one link
component to another when we encounter a marked vertex representing a
classical crossing involving those two link components. Each circuit of $C$
must cover all the link components that appear together in a connected
component of $\Gamma$, so if $\Gamma$ has $c(\Gamma)$ connected components
then $\mathcal{L}(D,C)$ must have at least $\mu-c(\Gamma)$ marked vertices. It
is a simple matter to construct an Euler circuit with precisely $\mu
-c(\Gamma)$ marked vertices: choose a spanning forest in $\Gamma$, and mark
one crossing of $D$ corresponding to each edge of the spanning forest.

\begin{lemma}
If $\mathcal{L}(D,C)$ has $\mu-c(\Gamma)$ marked vertices then no two of them
are adjacent to each other.
\end{lemma}

\begin{proof}
If two marked vertices $v$ and $w$ are neighbors then they are interlaced on a
circuit of $C$. If the circuit is $C_{1}=C_{11}vC_{12}wC_{13}vC_{14}w$ then
$C_{1}\ast v\ast w\ast v=$ $C_{11}vC_{14}wC_{13}vC_{12}$ is also an Euler
circuit of the same connected component of $U$. If we replace $C_{1}$ with
$C_{1}\ast v\ast w\ast v$ then we obtain an Euler system $C\ast v\ast w\ast v$
which has precisely the same marked vertices as $C$, except that $v$ and $w$
are no longer marked. This is impossible as $\mu-c(\Gamma)$ is the minimum
number of marked vertices.
\end{proof}

We are now ready to prove Theorem \ref{tangle}.

Suppose $D$ contains a tangle (more precisely, a \textit{2-tangle}), i.e., a
portion of $D$ that can be enclosed by a circle which intersects $D$ in
precisely four points, none of which is a crossing. Let $p_{1}$, $p_{2}$,
$p_{3}$ and $p_{4}$ be the four intersection points of $D$ and the tangle's
boundary circle. Let $C$ be a directed Euler system for $\vec{U}$, and
consider a circuit $C_{1}\in C$ that passes through $p_{1}$; suppose $C_{1}$
is directed into the circle at $p_{1}$. Clearly $C_{1}$ must leave the circle,
say at $p_{2}$. There are two possibilities: either $C_{1}$ also passes
through $p_{3}$ and $p_{4}$, or a different circuit of $C$ passes through them.

In the first case, suppose $C_{1}$ is $C_{11}p_{1}C_{12}p_{2}C_{13}p_{3}%
C_{14}p_{4}$, with $C_{11}$ and $C_{13}$ outside the circle. A vertex of
$\mathcal{L}(D,C)$ corresponds to a classical crossing of $D$; either the
crossing is inside the tangle, in which case the vertex cannot appear on
$C_{11}$ or $C_{13}$, or else the crossing is outside the tangle, in which
case the vertex cannot appear on $C_{12}$ or $C_{14}$. If the crossing appears
on more than one arc $C_{1j}$, then, it must appear either on $C_{11}$ and
$C_{13}$ or on $C_{12}$ and $C_{14}$. Clearly these two types of vertices are
interlaced with respect to $C_{1}$. Also, a vertex that appears only on a
single $C_{1j}$ cannot be interlaced with a vertex that does not appear on the
same $C_{1j}$. Consequently the first assertion of Theorem \ref{tangle} holds,
with the vertices that appear on both $C_{12}$ and $C_{14}$ adjacent to $a$ in
$H$ and the vertices that appear on both $C_{11}$ and $C_{13}$ adjacent to $a$
in $F$.

In the second case, no vertex corresponding to a crossing inside the tangle is
interlaced with any vertex corresponding to a crossing outside the tangle.
Consequently $\mathcal{L}(D,C)$ is the disjoint union of $F-a$ and $H-a$, so
it is a composition with $a$ an isolated vertex.

This completes the proof of the first assertion of Theorem \ref{tangle}. To
prove the second assertion, note that the lemma tells us that it is possible
to choose $C$ so that no two marked vertices of $\mathcal{L}(D,C)$ are
neighbors. As every neighbor of $a$ in $F$ is adjacent in $F\ast H$ to every
neighbor of $a$ in $H$, at least one of $F,H$ must contain no marked neighbor
of $a$.

\begin{corollary}
\label{tangle1} If $D$ contains a tangle then Theorem \ref{pjoin} may be used
to describe the Kauffman bracket of $D$ using appropriately weighted versions
of the subgraphs of $\mathcal{L}(D,C)$ induced by the vertices inside and
outside the tangle.
\end{corollary}

\begin{proof}
As discussed above, it is always possible to choose $C$ so that no two marked
vertices of $\mathcal{L}(D,C)$ are adjacent. Then Theorem \ref{tangle} tells
us that Theorem \ref{pjoin} applies, possibly with the names of $F$ and $H$ interchanged.
\end{proof}%

\begin{figure}
[h]
\begin{center}
\includegraphics[
trim=0.935096in 4.139700in 1.205565in 0.872867in,
height=4.3007in,
width=4.6086in
]%
{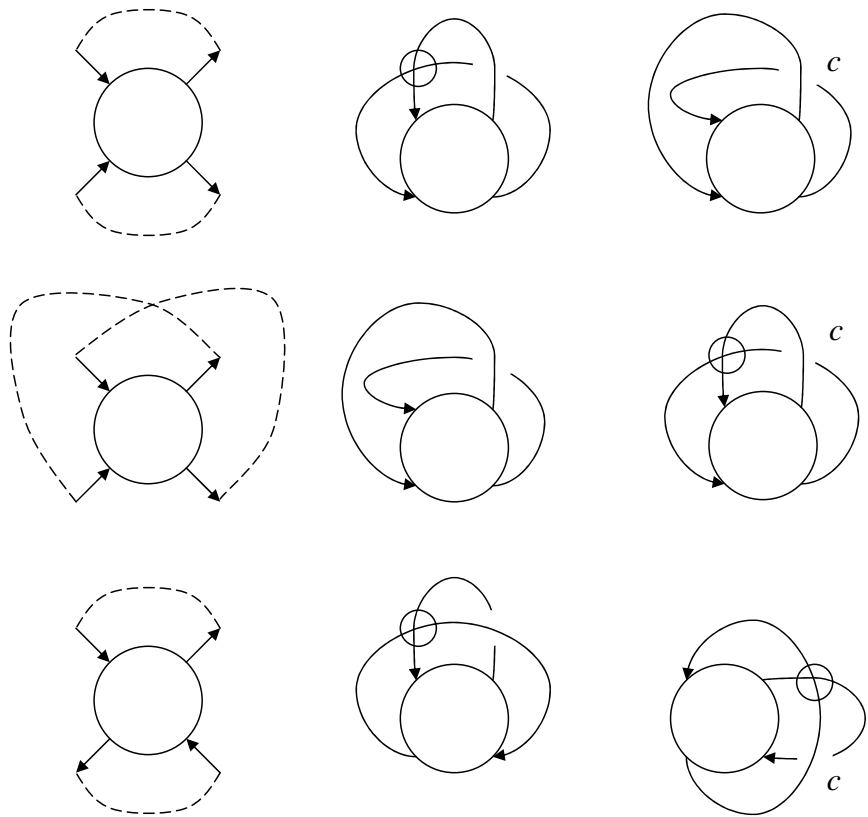}%
\caption{At left is a directed tangle in $D $, with dashes indicating the
Euler system $C$. $H-a$ is the subgraph of $\mathcal{L}(D,C)$ induced by
vertices corresponding to crossings inside the tangle. The center diagram
corresponds to $H^{\prime}$ and the right-hand diagram corresponds to
$H_{m}^{\prime}$.}%
\label{tangles4}%
\end{center}
\end{figure}

In the situation of Corollary \ref{tangle1}, Theorem \ref{pjoin} tells us that
$[D]$ is equal to the sum of the bracket polynomials of two graphs with
weighted vertices. It turns out that these two graphs are looped interlacement
graphs, so $[D]$ is actually equal to the sum of the bracket polynomials of
two link diagrams with weighted crossings. Suppose $D$ contains a directed
tangle as illustrated on the left-hand side of Figure \ref{tangles4}, and an
Euler system $C$ has been chosen so that the $F,H$ notations of Theorems
\ref{tangle} and \ref{pjoin} agree. (The dashed arcs in the figure indicate
the paths $C$ might follow as it leaves and re-enters the tangle.) Then $F-a$
is the subgraph of $\mathcal{L}(D,C)$ induced by the vertices corresponding to
the crossings outside the circle, $H-a$ is the subgraph of $\mathcal{L}(D,C)$
induced by the vertices corresponding to the crossings inside the circle, and
the neighbors of $a$ are the vertices of $F-a$ and $H-a$ that correspond to
crossings outside and inside the circle that are interlaced with each other.
No neighbor of $a$ in $H$ is marked. The looped interlacement graphs of the
diagrams shown in the center and on the right-hand side of each row of the
figure are then isomorphic to the graphs $H^{\prime}$ and $H_{m}^{\prime}$,
with the new classical crossings corresponding to $a$ and $a_{m}$
respectively. (The figure does not completely specify these two weighted
graphs, because the weights of $a$ and $a_{m}$ are not displayed.)

Observe that the portion of $D$ inside the circle is not disturbed in the two
new diagrams. Consequently if $D$ is obtained by substituting tangles for the
vertices of a 4-regular graph $P$ (i.e., $D=P\ast t_{1}...t_{\tau}$ in the
notation of \cite{C}), with $t_{\tau}$ containing the crossings outside the
circle, then $t_{1},...,t_{\tau-1}$ still appear as tangles in the two new
diagrams, and Theorem \ref{pjoin} may be applied to $t_{\tau-1},...,t_{1}$ in
turn. The result of applying Theorem \ref{pjoin} repeatedly is to express
$[D]$ as the sum of the bracket polynomials of $2^{\tau}$ crossing-weighted
$\tau$-crossing diagrams. Such a sum seems complicated but depending on the
structure of $D$, it may actually be considerably simpler than the definition
of $[D]$.%

\begin{figure}
[t]
\begin{center}
\includegraphics[
trim=1.061260in 6.407443in 1.080226in 0.811895in,
height=2.6429in,
width=4.6051in
]%
{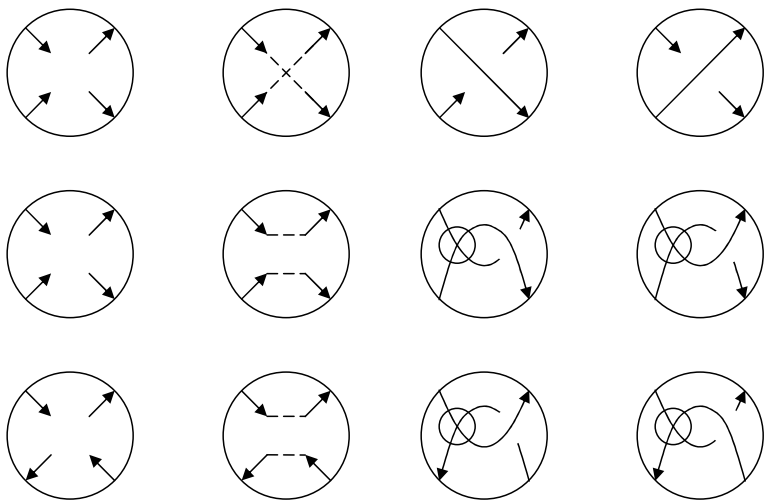}%
\caption{Each row represents a directed tangle located outside the circle, a
dashed Euler system in $D$, and associated diagrams $D_{+}$ and $D_{-}$.}%
\label{tangles3}%
\end{center}
\end{figure}

We close this section with the observation that if $F$ arises from a tangle in
a link diagram then it is possible to replace the graphs $F^{10}$ and $F^{01}$
of Corollary \ref{pjoin1} with graphs that also arise from link diagrams.
Suppose we are given a directed tangle as on the left-hand side of Figure
\ref{tangles3}. Note that the convention of the figure is the opposite of the
usual one -- the crossings of the tangle correspond to the vertices of the
subgraph $F-a$, so they are presumed to lie outside the circle; the four
segments inside the circle are the ends of arcs of the tangle. The dashes in
the second picture of each row indicate the path followed by $C$, and the
third and fourth pictures in that row indicate virtual diagrams we denote
$D_{+}$ and $D_{-}$ respectively. Let $F_{+}=\mathcal{L}(D_{+},C_{+})$ and
$F_{-}=\mathcal{L}(D_{-},C_{-})$, where $C_{+}$ and $C_{-}$ are the Euler
systems obtained from $C$ in the obvious ways. Also, let $v_{+}\in V(F_{+})$
and $v_{-}\in V(F_{-})$ be the vertices corresponding to the new classical
crossings; then $v_{-}$ is looped and $v_{+}$ is unlooped. If we weight
$v_{+}$ and $v_{-}$ with $\alpha(v_{+})=\alpha(v_{-})=A$ and $\beta
(v_{+})=\beta(v_{-})=B$, it follows that $[F_{+}]=A[F^{10}]+B[F^{01}]$ and
$[F_{-}]=B[F^{10}]+A[F^{01}]$. The formulas of Corollary \ref{pjoin1} then
imply the following.

\begin{corollary}
\label{pjoin3}
\[
(2-d-d^{2})\alpha(a)=-(d+1)[F]+\frac{[F_{+}]+[F_{-}]}{A+B}%
\]

\end{corollary}%

\[
(2-d-d^{2})\beta(a)=[F]+\frac{(A+B+Bd)[F_{+}]-(A+Ad+B)[F_{-}]}{A^{2}-B^{2}}%
\]

\[
(2-d-d^{2})\alpha(a_{m})=[F]+\frac{(A+B+Bd)[F_{-}]-(A+Ad+B)[F_{+}]}%
{A^{2}-B^{2}}%
\]

\section{Twin vertices}

Two edges of a graph incident on the same vertices are \textit{parallel}, and
two edges incident on a degree-two vertex are \textit{in series}. Series and
parallel edges in edge-weighted graphs can often be consolidated using some
appropriate combination of weights; for instance, in electrical circuit theory
two parallel resistors are equivalent to a single resistor with $R^{-1}%
=R_{1}^{-1}+R_{2}^{-1}$. Similarly, vertices $v$ and $w$ are called
\textit{twins} if they have the same neighbors outside $\{v,w\}$.%

\begin{figure}
[h]
\begin{center}
\includegraphics[
trim=1.606321in 8.842058in 1.607145in 0.805476in,
height=0.8155in,
width=3.8017in
]%
{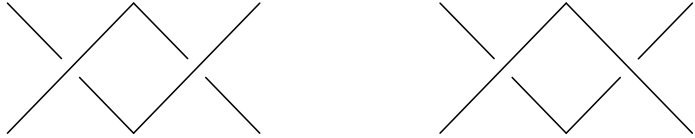}%
\caption{Crossings that are consecutive on both strands give rise to twin
vertices in $\mathcal{L}(D,C)$.}%
\label{twins}%
\end{center}
\end{figure}

\begin{figure}
[h]
\begin{center}
\includegraphics[
trim=1.605496in 8.432366in 1.610444in 0.671765in,
height=1.222in,
width=3.8017in
]%
{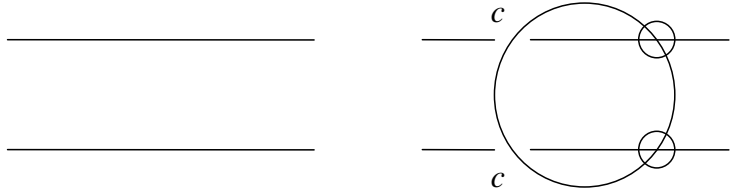}%
\caption{A construction that produces nonadjacent marked twins.}%
\label{twinex}%
\end{center}
\end{figure}

Some twin vertices in interlacement graphs arise from simple configurations in
link diagrams, like those pictured in Figure \ref{twins}. Other twin vertices
arise from more complicated configurations. For instance, let $D_{1} $ and
$D_{2}$ be two link diagrams with connected universes, each marked to identify
an Euler circuit. Suppose $D_{1}$ and $D_{2}$ are drawn together in the plane,
with some finite, positive number of transverse intersections. Let $D_{3}$ be
a link diagram obtained by considering each of these transverse intersections
to be a classical crossing. Then the universe $U_{3} $ of $D_{3}$ is
connected, but the marks on $D_{1}$ and $D_{2}$ do not describe an Euler
circuit in $U_{3}$. If we locate a portion of $D_{3}$ where an arc of $D_{1}$
lies parallel to an arc of $D_{2}$, and replace that portion of $D_{3}$ with
the configuration pictured on the right in Figure \ref{twinex}, then the
result is a marked link diagram $D$ whose marks do describe an Euler circuit
$C$. The two new classical crossings correspond to marked, nonadjacent twins
in $\mathcal{L}(D,C)$.%

In some situations, the weighted bracket polynomial allows for the
consolidation of twin vertices into one vertex.

\begin{proposition}
\label{lemma2} Let $v,w\in V(G)$ be unlooped twins.

(a) Suppose $v$ and $w$ are marked and not adjacent. Let $(G-w)^{\prime}$ be
the graph obtained from $G-w$ by changing the weights of $v$ to $\alpha
^{\prime}(v)=\alpha(v)\alpha(w)+\beta(v)\beta(w)$ and $\beta^{\prime
}(v)=\alpha(v)\beta(w)+\beta(v)\alpha(w)$. Then $[G]=[(G-w)^{\prime}].$

(b) Suppose $v$ and $w$ are marked and adjacent. Let $(G-w)^{\prime}$ be the
graph obtained from $G-w$ by changing the weights of $v$ to $\alpha^{\prime
}(v)=\alpha(v)\alpha(w)$ and $\beta^{\prime}(v)=\alpha(v)\beta(w)+\beta
(v)\alpha(w)+\beta(v)\beta(w)d$. Then $[G]=[(G-w)^{\prime}].$

(c) Suppose $v$ and $w$ are unmarked and adjacent. Let $(G-w)^{\prime}$ be the
graph obtained from $G-w$ by marking $v$ and changing the weights of $v$ to
$\alpha^{\prime}(v)=\alpha(v)\alpha(w)$ and $\beta^{\prime}(v)=\alpha
(v)\beta(w)+\beta(v)\alpha(w)+\beta(v)\beta(w)d$. Then $[G]=[(G-w)^{\prime}].$

(d) Suppose $v$ and $w$ are adjacent, with $v$ unmarked and $w$ marked. Let
$(G-w)^{\prime}$ be the graph obtained from $G-w$ by changing the weights of
$v$ to $\alpha^{\prime}(v)=\alpha(v)\alpha(w)$ and $\beta^{\prime}%
(v)=\alpha(v)\beta(w)+\beta(v)\alpha(w)+\beta(v)\beta(w)d$. Then
$[G]=[(G-w)^{\prime}].$
\end{proposition}

\begin{proof}
Observe that $v$ and $w$ give rise to nearly identical rows and columns of
$\mathcal{A}(G)$; they differ only in their common entries, and only if $v$
and $w$ are neighbors. The four parts of the lemma are all justified by
applying to $\mathcal{A}(G)_{T}$ the following nullity calculations over
$GF(2)$.%

\[
\nu%
\begin{pmatrix}
0 & 0 & \mathbf{1} & \mathbf{0}\\
0 & 0 & \mathbf{1} & \mathbf{0}\\
\mathbf{1} & \mathbf{1} & M_{11} & M_{12}\\
\mathbf{0} & \mathbf{0} & M_{21} & M_{22}%
\end{pmatrix}
-1=\nu%
\begin{pmatrix}
1 & 0 & \mathbf{1} & \mathbf{0}\\
0 & 0 & \mathbf{1} & \mathbf{0}\\
\mathbf{1} & \mathbf{1} & M_{11} & M_{12}\\
\mathbf{0} & \mathbf{0} & M_{21} & M_{22}%
\end{pmatrix}
=\nu%
\begin{pmatrix}
0 & \mathbf{1} & \mathbf{0}\\
\mathbf{1} & M_{11} & M_{12}\\
\mathbf{0} & M_{21} & M_{22}%
\end{pmatrix}
\]

\[
\nu%
\begin{pmatrix}
1 & 0 & \mathbf{1} & \mathbf{0}\\
0 & 1 & \mathbf{1} & \mathbf{0}\\
\mathbf{1} & \mathbf{1} & M_{11} & M_{12}\\
\mathbf{0} & \mathbf{0} & M_{21} & M_{22}%
\end{pmatrix}
=\nu%
\begin{pmatrix}
0 & 1 & \mathbf{1} & \mathbf{0}\\
1 & 0 & \mathbf{1} & \mathbf{0}\\
\mathbf{1} & \mathbf{1} & M_{11} & M_{12}\\
\mathbf{0} & \mathbf{0} & M_{21} & M_{22}%
\end{pmatrix}
=\nu%
\begin{pmatrix}
M_{11} & M_{12}\\
M_{21} & M_{22}%
\end{pmatrix}
\]

\[
\nu%
\begin{pmatrix}
1 & 1 & \mathbf{1} & \mathbf{0}\\
1 & 1 & \mathbf{1} & \mathbf{0}\\
\mathbf{1} & \mathbf{1} & M_{11} & M_{12}\\
\mathbf{0} & \mathbf{0} & M_{21} & M_{22}%
\end{pmatrix}
-1=\nu%
\begin{pmatrix}
1 & 1 & \mathbf{1} & \mathbf{0}\\
1 & 0 & \mathbf{1} & \mathbf{0}\\
\mathbf{1} & \mathbf{1} & M_{11} & M_{12}\\
\mathbf{0} & \mathbf{0} & M_{21} & M_{22}%
\end{pmatrix}
=\nu%
\begin{pmatrix}
1 & \mathbf{1} & \mathbf{0}\\
\mathbf{1} & M_{11} & M_{12}\\
\mathbf{0} & M_{21} & M_{22}%
\end{pmatrix}
\]

For instance, consider part (c). As neither $v$ nor $w$ is marked, the
definition of $\mathcal{A}(G)_{T}$ will not ever involve deleting the row and
column corresponding to either. The second nullity calculation shows that if
$T\subseteq V(G)$ has $v,w\notin T$ then the\ term in Definition
\ref{bracweigh}'s formula for $[G]$ corresponding to $T$ is the same as the
term in Definition \ref{bracweigh}'s formula for $[(G-w)^{\prime}]$
corresponding to $T$. The third nullity calculation shows that sum of
the\ terms corresponding to $T\cup\{v\}$, $T\cup\{w\}$ and $T\cup\{v,w\}$ in
Definition \ref{bracweigh}'s formula for $[G]$ coincides with the term in
Definition \ref{bracweigh}'s formula for $[(G-w)^{\prime}]$ corresponding to
$T\cup\{v\}$.
\end{proof}

\bigskip

Observe that the reduced graphs $(G-w)^{\prime}$ in parts (b) and (c) of
Proposition \ref{lemma2} are the same; this is not surprising, as the original
graphs differ only by a marked pivot. Also, parts (b)-(d) inductively imply
Theorem \ref{parallels} of the Introduction.

Two cases are missing from Proposition \ref{lemma2}. These cases are not fully
analogous to series-parallel reductions of edges, as they do not involve the
consolidation of two vertices into one vertex of one graph.

\begin{proposition}
\label{lemma3} Let $v,w\in V(G)$ be nonadjacent, unlooped twins.

(a) Suppose neither $v$ nor $w$ is marked. Let $(G-w)^{\prime}$ be the graph
obtained from $G-w$ by changing the weights of $v$ to $\alpha^{\prime
}(v)=\alpha(v)\alpha(w)d+\alpha(v)\beta(w)+\beta(v)\alpha(w)$ and
$\beta^{\prime}(v)=0$. Then $[G]=[(G-w)^{\prime}]+\beta(v)\beta(w)[G-v-w].$

(b) Suppose $v$ is unmarked and $w$ is marked. Let $(G-w)^{\prime}$ be the
graph obtained from $G-w$ by changing the weights of $v$ to $\alpha^{\prime
}(v)=\alpha(v)\alpha(w)+\alpha(v)\beta(w)$ and $\beta^{\prime}(v)=\beta
(v)\alpha(w)$. Then $[G]=[(G-w)^{\prime}]+\beta(v)\beta(w)[G-v-w].$
\end{proposition}

\begin{proof}
The proofs use the same nullity calculations that appear in the proof of
Proposition \ref{lemma2}. These cases are more complicated because each one
involves both the first nullity calculation and the second, and the
corresponding matrices do not appear together in $[G-v]$ or $[G-w]$.
\end{proof}

The inductive version of Proposition \ref{lemma3} allow us to simplify a graph
$G$ containing $k$ nonadjacent twins into one or two graphs with $\left\vert
V(G)\right\vert -k+1$ or $\left\vert V(G)\right\vert -k$ vertices. A
preliminary observation will be useful.

\begin{lemma}
\label{lemma1} Suppose two graphs $G_{1}$ and $G_{2}$ are identical except for
the weights of a single vertex $a$, and let $G$ be the graph that is identical
to both $G_{1}$ and $G_{2}$ except for $\alpha(a)=r_{1}\alpha_{1}%
(a)+r_{2}\alpha_{2}(a)$ and $\beta(a)=r_{1}\beta_{1}(a)+r_{2}\beta_{2}(a)$.
Then $[G]=r_{1}[G_{1}]+r_{2}[G_{2}]$.
\end{lemma}

\begin{proof}
This follows directly from Definition \ref{bracweigh}.
\end{proof}

\begin{corollary}
Let $v_{1},...,v_{k}$ be nonadjacent, unlooped twins in $G$. Then there are
weights $\alpha^{\prime}(v_{1}),\beta^{\prime}(v_{1})$ and a coefficient
$\gamma$ so that $[G]=[(G-v_{2}-...-v_{k})^{\prime}]+\gamma\cdot\lbrack
G-v_{1}-...-v_{k}]$.
\end{corollary}

\begin{proof}
Suppose first that $v_{1},...,v_{k}$ are all unmarked and $k\geq3$. Part (a)
of Proposition \ref{lemma3} gives us $[G]=[(G-v_{k})^{\prime}]+\beta
(v_{k-1})\beta(v_{k})[G-v_{k-1}-v_{k}]$. As $\beta^{\prime}(v_{k-1})=0$,
applying part (a) of Proposition \ref{lemma3} to $(G-v_{k})^{\prime}$ gives us
$[(G-v_{k})^{\prime}]=[((G-v_{k})^{\prime}-v_{k-1})^{\prime}]$, where the
weights of $v_{k-2}$ in $((G-v_{k})^{\prime}-v_{k-1})^{\prime}$ have been
changed to $\alpha^{\prime\prime}(v_{k-2})=\alpha^{\prime}(v_{k-1}%
)\alpha(v_{k-2})d+\beta(v_{k-2})\alpha^{\prime}(v_{k-1})$ and $\beta
^{\prime\prime}(v_{k-2})=0$. Lemma \ref{lemma1} tells us that $[G]=[((G-v_{k}%
)^{\prime}-v_{k-1})^{\prime}]+\beta(v_{k-1})\beta(v_{k})[G-v_{k-1}%
-v_{k}]=[(G-v_{k}-v_{k-1})^{\prime\prime\prime}]$, where the only difference
between $(G-v_{k}-v_{k-1})^{\prime\prime\prime}$ and $G-v_{k}-v_{k-1}$ is that
the weights of $v_{k-2}$ are different. Repeating this process $\left\lfloor
\frac{k-1}{2}\right\rfloor $ times will ultimately reduce $G$ to a version of
$G-v_{2}-...-v_{k}$ or $G-v_{3}-...-v_{k}$ which has been modified only by
changing the weights of $v_{1}$ or $v_{2}$. In the former case $\gamma=0$, and
in the latter case one more application of part (a) of Proposition
\ref{lemma3} is needed.

Suppose now that $k\geq3$ and $v_{k}$ is the only marked vertex among
$v_{1},...,v_{k}$. Part (b) of Proposition \ref{lemma3} gives us
$[G]=[(G-v_{k})^{\prime}]+\beta(v_{k-1})\beta(v_{k})[G-v_{k-1}-v_{k}]$.
Applying part (a) of Proposition \ref{lemma3} then tells us that $[G]$ is
\begin{gather*}
\beta^{\prime}(v_{k-1})\beta(v_{k-2})[(G-v_{k})^{\prime}-v_{k-1}-v_{k-2}]+\\
\lbrack((G-v_{k})^{\prime}-v_{k-1})^{\prime}]+\beta(v_{k-1})\beta
(v_{k})[G-v_{k-1}-v_{k}].
\end{gather*}
Lemma \ref{lemma1} tells us that the sum of the second and third terms is the
bracket polynomial of a graph that differs from $G-v_{k-1}-v_{k}$ only in the
weights of $v_{k-2}$. If $k>3$ then the paragraph above applies to this new
graph and also to $(G-v_{k})^{\prime}-v_{k-1}-v_{k-2}$, as neither has any
marked vertex among the remaining $v_{i}$. The resulting expressions may be
combined by using Lemma \ref{lemma1} to get appropriate weights for $v_{1}$
and simply adding together the coefficients multiplying $[G-v_{1}-...-v_{k}]$.

If there are two or more marked vertices among $v_{1},...,v_{k}$, then we may
apply part (a) of Proposition \ref{lemma2} repeatedly to bring the number of
marked vertices down to one.
\end{proof}

The first paragraph of the proof yields a partial \textquotedblleft
dual\textquotedblright\ of Theorem \ref{parallels}.

\begin{corollary}
\label{parallelsdual} Suppose $k\geq3$ is odd and $v_{1},...,v_{k}$ are
unmarked, nonadjacent, unlooped twins in $G$. Then $[G]=[(G-v_{2}%
-...-v_{k})^{\prime}]$ where
\[
\beta^{\prime}(v_{1})=\prod_{i=1}^{k}\beta(v_{i})\text{ \ and \ }%
\alpha^{\prime}(v_{1})=d^{-1}\left(  -\beta^{\prime}(v_{1})+\prod_{i=1}%
^{k}(\alpha(v_{i})d+\beta(v_{i}))\right)  .
\]

\end{corollary}

\section{A recursion for the weighted bracket}

A recursion for the bracket polynomial of marked graphs was given in
\cite{T1}. Modifying this recursion to describe the vertex-weighted version of
the bracket is not difficult. Recall that the \textit{open neighborhood}
$N(v)$ of a vertex of a graph contains the vertices $w\neq v$ such that
$\{v,w\}\in E(G)$.

\begin{definition}
If $v\in V(G)$ then the \emph{local complement} $G^{v}$ has $V(G^{v})=V(G)$
and $E(G^{v})=\{\{a,b\}$
$\vert$
$a\notin N(v)$ and $\{a,b\}\in E(G)\}$ $\cup$ $\{\{a,b\}$
$\vert$
$a,b\in N(v)$ and $\{a,b\}\notin E(G)\}$.
\end{definition}

That is, $G^{v}$ is obtained from $G$ by toggling loops and non-loop edges in
the subgraph induced by $N(v)$.

\begin{definition}
If $v,w\in V(G)$ then the \emph{pivot} $G^{vw}$ is obtained from $G$ by
toggling every adjacency between vertices $a$ and $b$ such that $a\in N(v)$,
$b\in N(w)$, and either $a\notin N(w)$ or $b\notin N(v)$.
\end{definition}

\begin{definition}
If $v,w\in V(G)$ are neighbors then the \emph{marked pivot} $G_{c}^{vw}$ is
obtained from $G^{vw}$ by interchanging the neighbors of $v$ and $w$ and
toggling the marks on $v$ and $w$.
\end{definition}

None of these three operations affects free loops or the weights of any vertex.

\begin{theorem}
\label{recursion} The weighted bracket polynomial of $G$ satisfies the following.

(a) Suppose $G$ has $\phi$ free loops, and let $G^{\prime}$ be the graph
obtained from $G$ by removing the free loops. Then $[G]=d^{\phi}\cdot\lbrack
G^{\prime}]$.

(b) Suppose $v$ is a looped vertex of $G$, and let $G^{\prime}$ be the graph
obtained from $G$ by removing the loop on $v$ and interchanging $\alpha(v)$
and $\beta(v)$. Then $[G]=[G^{\prime}]$.

(c) If $v$ and $w$ are neighbors in $G$ then $[G]=[G_{c}^{vw}]$.

(d) Suppose $v$ is unlooped and marked, and no neighbor of $v$ is marked.
Then
\[
\lbrack G]=\alpha(v)[G-v]+\beta(v)[G^{v}-v],
\]
where $G-v$ is obtained from $G$ by removing $v$ and every edge incident on
$v$.

(e) Let $v$ and $w$ be adjacent, unlooped, unmarked vertices. If no neighbor
of $v$ is marked then
\[
\lbrack G]=\alpha(v)\alpha(w)[G^{vw}-v-w]+\alpha(v)\beta(w)[(G^{vw}%
)^{v}-v-w]+\beta(v)[G^{v}-v].
\]

(f) If $v$ is an isolated, unlooped, unmarked vertex of $G$ then
\[
\lbrack G]=(\alpha(v)d+\beta(v))\cdot\lbrack G-v]\text{.}%
\]

(g) The empty graph $\emptyset$ has $[\emptyset]=1$.
\end{theorem}

\begin{proof}
Parts (a), (b), (f) and (g) follow immediately from Definition \ref{bracweigh}%
. Part (c) follows from the fact that for every $T\subseteq V(G)$, the
$GF(2)$-nullities of $\mathcal{A}(G)_{T}$ and $\mathcal{A}(G_{c}^{vw})_{T}$
are the same; see Section 5 of \cite{T1}. Parts (d) and (e) are proven just as
in the unweighted case; see Section 6 of \cite{T1}.
\end{proof}

\bigskip

Theorem \ref{recursion} provides a recursive algorithm for the weighted
bracket polynomial: first use (a) to eliminate free loops; then use (b) to
eliminate loops; then use (c) to eliminate adjacencies between marked
vertices; then use (d) to eliminate marked vertices; then use (e) to eliminate
the remaining adjacencies; and finally use (f) and (g) to calculate the
bracket polynomials of the remaining edgeless, unmarked graphs. Different
individual implementations of the algorithm will involve applying parts
(c)-(e) at different locations in $G$, and just as in the unweighted case
there is no canonical way to find the most efficient implementation.

Another property of the unweighted bracket that extends directly to the
weighted version is this: if $G$ is the union of disjoint subgraphs $G_{1}$
and $G_{2}$ then $[G]=[G_{1}]\cdot\lbrack G_{2}]$.

\section{Proof of Theorem \ref{pjoin}}

In outline, our proof of Theorem \ref{pjoin} is similar to the corresponding
proof for vertex-weighted interlace polynomials \cite{Tw}. This similarity is
not surprising as weighted interlace polynomials can also be calculated
recursively using pivots and local complements. A more complicated argument is
required here, however, because interlace polynomials have more convenient
reductions for twin and pendant vertices, and do not involve marked vertices.

Suppose that $G=F\ast H$, and $a$ has no marked neighbor in $H$. Consider the
following calculation of $[G]$.

First, unloop each looped vertex $v\in V(F)$, and reverse $\alpha(v)$ and
$\beta(v)$. This does not affect any edges other than loops in $F-a$, so when
we are done we have a composition $F^{\prime}\ast H$ with no loops in
$F^{\prime}$.

Second, eliminate adjacencies between marked vertices of $F^{\prime}-a$ using
marked pivots. As no two vertices of $H$ have different nonempty sets of
neighbors in $F^{\prime}$, these marked pivots will not affect the internal
structure of $H$. However there may be extensive changes within $F^{\prime}$.
The graph $F^{\prime\prime}\ast H$ resulting from this step is not unique; but
any two differ only by marked pivots, and have the same bracket polynomial as
$G$.

As no neighbor of $a$ in $H$ is marked, no marked vertex of $F^{\prime\prime}$
has a marked neighbor in $F^{\prime\prime}\ast H$. Consequently part (d) of
Theorem \ref{recursion} may be used to eliminate every marked vertex in
$F^{\prime\prime}$. The result is a sum of terms, each of which is the product
of an initial factor and some bracket polynomial $[F^{\prime\prime\prime}\ast
H]$ or $[F^{\prime\prime\prime}\ast H^{a}]$.

Fourth, eliminate all remaining adjacencies in each $F^{\prime\prime\prime}-a$
using part (e) of Theorem \ref{recursion}. As each resulting term is obtained
by applying local complement or pivot operations, each term is the product of
an initial factor and a bracket polynomial $[F^{\prime\prime\prime\prime}\ast
H]$ or $[F^{\prime\prime\prime\prime}\ast H^{a}]$, with $F^{\prime\prime
\prime\prime}$ unmarked and edgeless. Consequently every vertex not adjacent
to $a$ in any $F^{\prime\prime\prime\prime}$ is now unmarked and isolated;
part (f) of Theorem \ref{recursion} tells us that each term is unchanged if we
remove all such vertices from that term's $F^{\prime\prime\prime\prime}$ and
multiply the initial factor appropriately. If a term $[F^{\prime\prime
\prime\prime}\ast H]$ involves a graph $F^{\prime\prime\prime\prime}$ in which
no vertex neighbors $a$, we may replace that term with a bracket polynomial
$[P_{1m}\ast H]$, where $V(P_{1m})=\{a,a_{m}\}$, $a$ and $a_{m}$ are adjacent,
and $a_{m}$ is a new marked vertex with $\beta(a_{m})=0$ and $\alpha(a_{m}%
)=1$. Similarly, any term $[F^{\prime\prime\prime\prime}\ast H^{a}]$ in which
$F^{\prime\prime\prime\prime}$ contains no neighbor of $a$ may be replaced
with a bracket polynomial $[\tilde{P}_{1m}\ast H^{a}]$, where $V(\tilde
{P}_{1m})=\{a,\tilde{a}_{m}\}$, $a$ and $\tilde{a}_{m}$ are adjacent, and
$\tilde{a}_{m}$ is a new marked vertex with $\beta(\tilde{a}_{m})=0$ and
$\alpha(\tilde{a}_{m})=1$. After these manipulations $[G]$ is expressed as a
sum of terms each of which is the product of an initial factor and a weighted
bracket polynomial $[F^{\prime\prime\prime\prime}\ast H]$ or $[F^{\prime
\prime\prime\prime}\ast H^{a}]$ in which $F^{\prime\prime\prime\prime}-a$ is a
nonempty collection of isolated unlooped neighbors of $a$.

Fifth, apply Proposition \ref{lemma3} and part (a) of Proposition \ref{lemma2}
repeatedly\ to each term in which $F^{\prime\prime\prime\prime}-a$ contains
more than one vertex, ultimately obtaining an expression of $[G]$ as a sum of
terms each of which is the product of an initial factor with a weighted
bracket $[F^{!}\ast H]$ or $[F^{!}\ast H^{a}]$ where $F^{!}-a $ consists of a
single unlooped neighbor of $a$, denoted $a_{u}$, $a_{m}$, $\tilde{a}_{u}$ or
$\tilde{a}_{m}$ according to whether that vertex is unmarked or marked
(subindex $u$ or $m$) and whether that term involves $H $ or $H^{a}$ (the
latter indicated by tilde). The value of an individual term is not changed if
we multiply both the $\alpha$ and $\beta$ weights of that term's vertex
$a_{x}$ or $\tilde{a}_{x}$ by the initial factor, and then replace the initial
factor with 1; consequently we may as well presume that the initial factors
are all 1. Using Lemma \ref{lemma1}, we see that if we sum the $\alpha$ and
$\beta$ weights of $a_{u}$, $\tilde{a}_{u}$, $a_{m} $ and $\tilde{a}_{m}$ in
all the different terms involving each, then $[G]$ is equal to the total of
four individual weighted bracket polynomials:%

\begin{equation}
\lbrack G]=[P_{1}\ast H]+[\tilde{P}_{1}\ast H^{a}]+[P_{1m}\ast H]+[\tilde
{P}_{1m}\ast H^{a}].\label{joine}%
\end{equation}
Here each of the graphs $P_{1},\tilde{P}_{1},P_{1m},\tilde{P}_{1m}$ consists
of $a$ and a single neighbor of $a$, denoted $a_{u}$, $\tilde{a}_{u}$, $a_{m}
$, $\tilde{a}_{m}$ respectively as before.

We now seem to have the job of determining eight unknowns, namely the $\alpha$
and $\beta$ weights of the four vertices $a_{u}$, $\tilde{a}_{u}$, $a_{m}$,
$\tilde{a}_{m}$. It turns out though that five of these unknowns are not necessary.

For instance, $\beta(\tilde{a}_{m})$ appears in terms of $[\tilde{P}_{1m}\ast
H^{a}]$ that involve the $GF(2)$-nullities of matrices of the form
\[%
\begin{pmatrix}
1 & \mathbf{0} & \mathbf{1}\\
\mathbf{0} & M_{11} & M_{12}\\
\mathbf{1} & M_{21} & M_{22}%
\end{pmatrix}
,\mathrm{~where~}
\begin{pmatrix}
M_{11} & M_{12}\\
M_{21} & M_{22}%
\end{pmatrix}
=\mathcal{A}(H^{a})_{T}.
\]
As $H$ contains no marked neighbor of $a$, for each such matrix
\[
\nu%
\begin{pmatrix}
1 & \mathbf{0} & \mathbf{1}\\
\mathbf{0} & M_{11} & M_{12}\\
\mathbf{1} & M_{21} & M_{22}%
\end{pmatrix}
=\nu%
\begin{pmatrix}
M_{11} & M_{12}\\
M_{21} & \bar{M}_{22}%
\end{pmatrix}
=\nu(\mathcal{A}(H)_{T}),
\]
where the overbar indicates that every entry of $\bar{M}_{22}$ is different
from the corresponding entry of $M_{22}$. (N.b. If we did not know that the
neighbors of $a$ in $H$ are unmarked, the last equality would be suspect as
the definition of $\mathcal{A}(H)_{T}$ would involve removing some rows of
$\bar{M}_{22}$ and retaining some rows removed in the definition of
$\mathcal{A}(H^{a})_{T}$.) Consequently the contribution to the sum of
(\ref{joine}) made by the terms in which $\beta(\tilde{a}_{m})$ appears may be
provided equally well by terms of $[P_{1m}\ast H]$ in which $\alpha(a_{m})$
appears. That is, the sum is unchanged if we replace $\alpha(a_{m})$ by
$\alpha(a_{m})+\beta(\tilde{a}_{m})$ and replace $\beta(\tilde{a}_{m})$ by $0$.

The terms of $[\tilde{P}_{1}\ast H^{a}]$ in which $\beta(\tilde{a}_{u})$
appears involve the $GF(2)$-nullities of precisely the same matrices just
discussed. Consequently the contributions to the sum of (\ref{joine}) made by
the terms of $[\tilde{P}_{1}\ast H^{a}]$ in which $\beta(\tilde{a}_{u})$
appears may also be provided by the terms of $[P_{1m}\ast H]$ in which
$\alpha(a_{m})$ appears; that is, we may replace $\alpha(a_{m})$ by
$\beta(\tilde{a}_{u})+\alpha(a_{m})$ and replace $\beta(\tilde{a}_{u})$ by $0$.

The equality
\[
\nu%
\begin{pmatrix}
M_{11} & M_{12}\\
M_{21} & M_{22}%
\end{pmatrix}
=\nu%
\begin{pmatrix}
1 & \mathbf{0} & \mathbf{1}\\
\mathbf{0} & M_{11} & M_{12}\\
\mathbf{1} & M_{21} & \bar{M}_{22}%
\end{pmatrix}
\]
tells us that the contributions of the terms in which $\alpha(\tilde{a}_{m}) $
appears can be duplicated by the terms in which $\beta(a_{u})$ appears, i.e.,
the sum of (\ref{joine}) is unchanged if we replace $\beta(a_{u})$ by
$\beta(a_{u})+\alpha(\tilde{a}_{m})$ and replace $\alpha(\tilde{a}_{m})$ by
$0$. With both weights of $\tilde{a}_{m}$ now $0$, $[\tilde{P}_{1m}\ast
H^{a}]=0$ makes no contribution to (\ref{joine}).

In the same way, the equality
\[
\nu%
\begin{pmatrix}
0 & \mathbf{0} & \mathbf{1}\\
\mathbf{0} & M_{11} & M_{12}\\
\mathbf{1} & M_{21} & \bar{M}_{22}%
\end{pmatrix}
=\nu%
\begin{pmatrix}
0 & \mathbf{0} & \mathbf{1}\\
\mathbf{0} & M_{11} & M_{12}\\
\mathbf{1} & M_{21} & M_{22}%
\end{pmatrix}
\]
tell us that $\tilde{a}_{u}$ is not needed: the contributions of the terms of
(\ref{joine}) involving $\alpha(\tilde{a}_{u})$ may be provided by the terms
involving $\alpha(a_{u})$.

The terms of (\ref{joine}) in which $\beta(a_{m})$ appears involve the same
nullities
\[
\nu%
\begin{pmatrix}
1 & \mathbf{0} & \mathbf{1}\\
\mathbf{0} & M_{11} & M_{12}\\
\mathbf{1} & M_{21} & M_{22}%
\end{pmatrix}
\]
that appear in the terms in which $\beta(a_{u})$ appears, so the sum of
(\ref{joine}) is unchanged if we replace $\beta(a_{u})$ with $\beta
(a_{u})+\beta(a_{m})$ and replace $\beta(a_{m})$ with $0$.

$P_{1}\ast H$ and $P_{1m}\ast H$ are isomorphic to $H^{\prime}$ and
$H_{m}^{\prime}$ respectively, so the statement of Theorem \ref{pjoin} follows.

\section{Subset formulas for the weights of Theorem \ref{pjoin}}

In this section we use linear algebra over $GF(2)$ to derive formulas for the
weights $\alpha(a_{u})$, $\beta(a_{u})$ and $\alpha(a_{m})$ of Theorem
\ref{pjoin} from Definition \ref{bracweigh}.

Suppose $T\subseteq V(F-a)$, and let $i_{1},...,i_{k}$ be the indices of the
rows and columns not removed from $\mathcal{A}(F-a)$ in obtaining
$\mathcal{A}(F-a)_{T}$. That is, $V(F-a)$ -- $\{v_{i_{1}},...,v_{i_{k}}%
\}=\{$marked $v\in V(F-a)|$ either $v\in T$ is looped or $v\notin T$ is not
looped$\}$. Let $\rho=(\rho_{i_{1}}...\rho_{i_{k}})$ be the row vector with
$\rho_{i_{j}}=1$ if $v_{i_{j}}$ is a neighbor of $a$, and let $\kappa$ be the
column vector obtained by transposing $\rho$. According to Lemma 2 of
\cite{BBCP}, the three nullities
\[
\nu\left(
\begin{array}
[c]{cc}%
\mathcal{A}(F)_{T} & \kappa\\
\rho & 0
\end{array}
\right)  ,\nu\left(
\begin{array}
[c]{cc}%
\mathcal{A}(F)_{T} & \kappa\\
\rho & 1
\end{array}
\right)  ,\nu\left(  \mathcal{A}(F)_{T}\right)
\]
are of the form $\nu+1$, $\nu$, $\nu$ in some order. We say $T$ is of
\emph{type} 1, 2 or 3 according to whether the nullity $\nu+1$ appears first,
second or third.

\begin{figure}
[ptb]
\begin{center}
\includegraphics[
trim=1.335027in 8.567147in 0.537639in 1.072899in,
height=0.8198in,
width=4.8084in
]%
{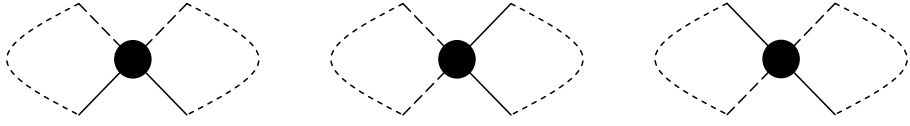}%
\caption{The circuit partition pictured in the middle contains one more
circuit.}%
\label{nullity}%
\end{center}
\end{figure}

If $D$ is a link diagram, then the three displayed matrices correspond to
three circuit partitions in $U$. At every vertex other than $a$, the three
partitions have the same transition. In this situation Lemma 2 of \cite{BBCP}
asserts that two of the partitions contain the same number of circuits, and
the third contains one more circuit. An example of type 2 is illustrated
schematically in Figure \ref{nullity}.%

For each subset $T\subseteq V(F-a)$ let $contr(T)$ denote the contribution of
$T$ to Definition \ref{bracweigh}'s formula for $[F-a]$,
\[
(\prod_{v\notin T}\alpha(v))(\prod_{t\in T}\beta(t))d^{\nu(\mathcal{A}%
(G)_{T})}.
\]

If $T$ is of type 1 then the sum of the contributions of $T$ and
$T\cup\{v_{10}\}$ to $[F^{10}]$ is $d\cdot contr(T)$, while the sum of the
contributions of $T$ and $T\cup\{v_{01}\}$ to $[F^{01}]$ is $contr(T)$. If $T$
is of type 2 then the sum of the contributions of $T$ and $T\cup\{v_{10}\}$ to
$[F^{10}]$ is $contr(T)$, and the sum of the contributions of $T$ and
$T\cup\{v_{01}\}$ to $[F^{01}]$ is $d\cdot contr(T)$. If $T$ is of type 3 then
the sum of the contributions of $T$ and $T\cup\{v_{10}\}$ to $[F^{10}]$ is
$(1/d)\cdot contr(T)$, and the sum of the contributions of $T$ and
$T\cup\{v_{01}\}$ to $[F^{01}]$ is also $(1/d)\cdot contr(T)$. As noted in the
Introduction, Theorem \ref{pjoin} implies these equations.
\begin{align*}
(2-d-d^{2})\alpha(a)  & =-(d+1)[F-a]+[F^{10}]+[F^{01}]\\
(2-d-d^{2})\beta(a)  & =[F-a]+[F^{10}]-(d+1)[F^{01}]\\
(2-d-d^{2})\alpha(a_{m})  & =[F-a]-(d+1)[F^{10}]+[F^{01}]
\end{align*}
It follows that if we denote the total of the contributions of the sets of
type $i$ to $[F-a]$ by $contr_{i}$, then these equalities hold.
\begin{align*}
(2-d-d^{2})\alpha(a)  & =(-(d+1)+(2/d))\cdot contr_{3}\\
(2-d-d^{2})\beta(a)  & =(2-d(d+1))\cdot contr_{2}\\
(2-d-d^{2})\alpha(a_{m})  & =(2-d(d+1))\cdot contr_{1}%
\end{align*}

\begin{corollary}
\label{pjoin4} The weights of Theorem \ref{pjoin} are $\alpha(a_{m}%
)=contr_{1}$, $\beta(a)=contr_{2}$ and $\alpha(a)=(contr_{3})/d$, where
\[
contr_{i}=\sum_{\substack{T\subseteq V(F-a) \\\mathrm{of\,\,type\,\,}i}%
}(\prod_{v\notin T}\alpha(v))(\prod_{t\in T}\beta(t))d^{\nu(\mathcal{A}%
(G)_{T})}.
\]

\end{corollary}

\section{Some examples}

The simplest example of Theorem \ref{pjoin} involves the two-vertex graph $F$
in which $a$ and $v$ are unlooped, unmarked neighbors. Let $\alpha=\alpha(v)$
and $\beta=\beta(v)$. For this graph Corollary \ref{pjoin1} gives the
following.
\begin{align*}
(2-d-d^{2})\alpha(a)  & =-(d+1)(\alpha d+\beta)+\alpha+\beta+\alpha+\beta
d=(2-(d+1)d)\alpha\\
(2-d-d^{2})\beta(a)  & =\alpha d+\beta+\alpha+\beta-(d+1)(\alpha+\beta
d)=(2-(d+1)d)\beta\\
(2-d-d^{2})\alpha(a_{m})  & =\alpha d+\beta-(d+1)(\alpha+\beta)+\alpha+\beta
d=0
\end{align*}
The assertion of Theorem 1 is then trivial, as $F\ast H$ and $H^{\prime}$ are
identical graphs.

A slightly more complicated example involves the same graph $F$, but with $v$
marked; Corollary \ref{pjoin1} gives the following.
\begin{align*}
(2-d-d^{2})\alpha(a)  & =-(d+1)(\alpha+\beta)+\alpha d+\beta+\alpha+\beta
d=0\\
(2-d-d^{2})\beta(a)  & =\alpha+\beta+\alpha d+\beta-(d+1)(\alpha+\beta
d)=(2-(d+1)d)\beta\\
(2-d-d^{2})\alpha(a_{m})  & =\alpha+\beta-(d+1)(\alpha d+\beta)+\alpha+\beta
d=(2-(d+1)d)\alpha
\end{align*}
This seems more interesting than the result of the first example, but it is
just as trivial. In the last step of the proof of Theorem \ref{pjoin} we see
that setting $\beta(a_{m})=0$ as in the statement of the theorem is arbitrary;
any choice of $\beta(a_{m})$ will satisfy the theorem, so long as the sum
$\beta(a)+\beta(a_{m})$ is correct. With $\beta(a)=0$ instead, the assertion
of Theorem \ref{pjoin} simply acknowledges that in this example, $F\ast H$ and
$H_{m}^{\prime}$ are identical.%

\begin{figure}
[ptbh]
\begin{center}
\includegraphics[
trim=0.667926in 8.768250in 1.338326in 0.676044in,
height=0.9669in,
width=4.708in
]%
{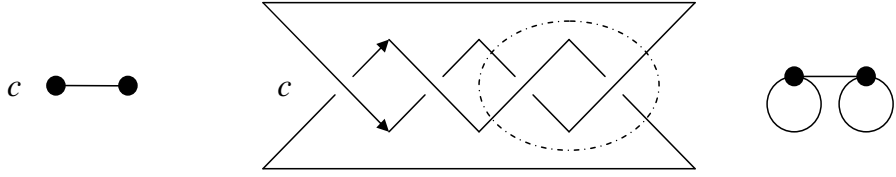}%
\caption{A two-component unlink. $F-a$ is on the left and $H-a$ is on the
right.}%
\label{example1}%
\end{center}
\end{figure}

Let $D$ be the link diagram pictured in the middle of Figure \ref{example1}.
Then $\mathcal{L}(D,C)=F\ast H$, where $F-a,H-a$ are the pictured two-vertex
graphs and $a$ is adjacent to all four of these vertices. Consider both $F-a$
and $H-a$ to have the standard weight functions $\alpha\equiv A$ and
$\beta\equiv B$. Corollary \ref{pjoin1} tells us that $\alpha(a)=A^{2}$,
$\beta(a)=2AB+B^{2}d$ and $\alpha(a_{m})=0$. Theorem \ref{pjoin} then asserts
that
\begin{gather*}
\lbrack F\ast H]=[H^{\prime}]=\\
A^{2}\cdot(A^{2}d+2AB+B^{2}d)+(2AB+B^{2}d)\cdot(A^{2}d^{2}+2ABd+B^{2}).
\end{gather*}
(The same conclusion follows from Theorem \ref{parallels}.) The writhe of $D$
is $0$, so this bracket polynomial determines the Jones polynomial through the
following two-stage evaluation. First $d\mapsto-A^{2}-B^{2}$ and $B\mapsto
A^{-1}$ yield
\[
f_{D}(A)=A^{2}\cdot(-A^{4}-A^{-4})+(1-A^{-4})\cdot(A^{6})=-A^{-2}-A^{2},
\]
and then $V=f_{D}(t^{-1/4})=-t^{1/2}-t^{-1/2}$, as we would expect for the
two-component unlink.

Suppose we modify the diagram in Figure \ref{example1} by reversing the
crossing on the right, effectively removing one loop from $H-a$. Then
\begin{gather*}
\lbrack F\ast H]=[H^{\prime}]=\\
A^{2}\cdot(A^{2}+2ABd+B^{2})+(2AB+B^{2}d)\cdot(A^{2}d+AB+ABd^{2}+B^{2}d).
\end{gather*}
The writhe of the diagram is now 2, so the Jones polynomial is obtained by
first calculating $f_{D}(A)$:
\begin{gather*}
A^{-6}\cdot(A^{4}+2(-A^{4}-1)+1+(1-A^{-4})\cdot(-A^{4}+A^{4}+2+A^{-4}%
-1-A^{-4}))\\
=A^{-6}\cdot(-A^{4}-A^{-4}),
\end{gather*}
and then evaluating $f_{D}(t^{-1/4})=t^{3/2}\cdot(-t-t^{-1})$, the correct
value for the positive Hopf link.

\begin{figure}
[t]
\begin{center}
\includegraphics[
trim=0.935921in 8.695510in 0.937570in 0.672835in,
height=1.0239in,
width=4.8075in
]%
{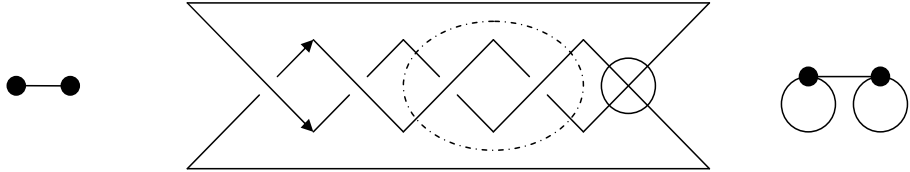}%
\caption{An unknot. $F-a$ is on the left and $H-a$ is on the right.}%
\label{example3}%
\end{center}
\end{figure}

Suppose instead we modify the graph \thinspace$F-a$\ pictured in Figure
\ref{example1} by removing the mark. Corollary \ref{pjoin1} tells us that
$\alpha(a)=0$, $\beta(a)=2AB+B^{2}d$ and $\alpha(a_{m})=A^{2}$. According to
Theorem \ref{pjoin},
\begin{gather*}
\lbrack F\ast H]=[H^{\prime}]+[H_{m}^{\prime}]=\\
(2AB+B^{2}d)\cdot(A^{2}d^{2}+2ABd+B^{2})+A^{2}\cdot(A^{2}d+2AB+B^{2}).
\end{gather*}
As the writhe is 0, this yields
\[
f_{D}(A)=(1-A^{-4})\cdot(A^{6})+A^{2}\cdot(1-A^{4}+A^{-2})=1.
\]
This seems incorrect at first glance, because Figure \ref{example1} displays a
two-component link. Note however that although the graph obtained by removing
the mark from $F-a$ is certainly a legitimate marked graph, we cannot
legitimately remove the mark from the link diagram in Figure \ref{example1},
because the marks on a connected link diagram must identify an Euler circuit.
Figure \ref{example3} exhibits an unmarked link diagram with the appropriate
graphs $F-a$ and $H-a$; it is an unknot rather than a two-component unlink, so
$V=1$ is indeed its Jones polynomial.%

\begin{figure}
[h]
\begin{center}
\includegraphics[
trim=0.668751in 8.766110in 0.937570in 0.669626in,
height=0.9738in,
width=5.0081in
]%
{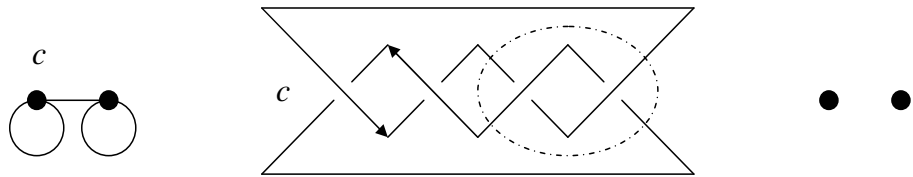}%
\caption{A two-component unlink. The marked vertex on the left is the only
neighbor of $a$ in $F$; $a$ is also adjacent to the two vertices of $H-a$ on
the right.}%
\label{example2}%
\end{center}
\end{figure}

Suppose we now modify the diagram in Figure \ref{example1} by reversing the
orientation of one component, as in Figure \ref{example2}. The only vertex of
$F$ adjacent to $a$ is the marked vertex, while both vertices of $H-a$ are
neighbors of $a$. Corollary \ref{pjoin1} tells us that $\alpha(a)=A^{2}$,
$\beta(a)=0$ and $\alpha(a_{m})=2AB+B^{2}d$. Then Theorem \ref{pjoin} asserts
that
\begin{gather*}
\lbrack F\ast H]=[H^{\prime}]+[H_{m}^{\prime}]=\\
A^{2}\cdot(A^{2}d+2AB+B^{2}d)+(2AB+B^{2}d)\cdot(A^{2}d^{2}+2ABd+B^{2}),
\end{gather*}
just as in the discussion of Figure \ref{example1} above. This is not
surprising, as the Kauffman bracket of a link diagram is independent of the
orientations of the link components.%

In our last example we consider the unmarked version of the example pictured
in Figure \ref{tanglex}. Direct calculations yield the following.%

\begin{figure}
[ptb]
\begin{center}
\includegraphics[
trim=3.077407in 7.893243in 2.275072in 0.406482in,
height=1.8256in,
width=2.1975in
]%
{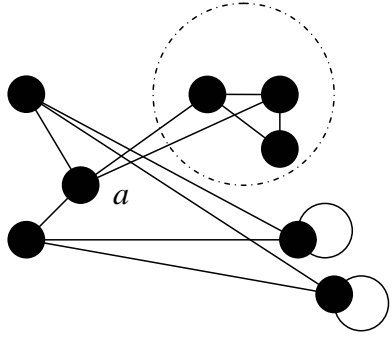}%
\caption{The unmarked graphs $F$ and $H$ from Figure \ref{tanglex}, with $H-a$
inside the dashed circle and $F $ outside the dashed circle.}%
\label{tanglex2}%
\end{center}
\end{figure}
\begin{align*}
\lbrack F-a]  & =A^{4}d^{2}+4A^{3}Bd+5A^{2}B^{2}+A^{2}B^{2}d^{2}%
+4AB^{3}d+B^{4}d^{2}\\
\lbrack F^{01}]  & =A^{4}d+2A^{3}B+2A^{3}Bd+5A^{2}B^{2}+A^{2}B^{2}%
d^{2}+4AB^{3}d+B^{4}d^{2}\\
\lbrack F^{10}]  & =A^{4}d+2A^{3}B+2A^{3}Bd^{2}+5A^{2}B^{2}d+A^{2}B^{2}%
d^{3}+4AB^{3}d^{2}+B^{4}d^{3}\\
\lbrack H-a]  & =A^{3}d+3A^{2}B+3AB^{2}d+B^{3}d^{2}\\
\lbrack H^{01}]  & =A^{3}d+3A^{2}B+AB^{2}d+2AB^{2}+B^{3}d\\
\lbrack H^{10}]  & =A^{3}d^{2}+3A^{2}Bd+AB^{2}d^{2}+2AB^{2}+B^{3}d
\end{align*}
Corollary \ref{pjoin1} gives us these weights.
\begin{align*}
\alpha(a(F))  & =A^{4}d+2A^{3}B,\beta(a(F))=0,\\
\alpha(a_{m}(F))  & =2A^{3}Bd+A^{2}B^{2}(5+d^{2})+4AB^{3}d+B^{4}d^{2},\\
\alpha(a(H))  & =2AB^{2}+B^{3}d,\beta(a(H))=0,\alpha(a_{m}(H))=A^{3}%
d+3A^{2}B+AB^{2}d
\end{align*}
Corollary \ref{pjoin2} tells us that
\begin{gather*}
\lbrack F\ast H]=(A^{4}d+2A^{3}B)\cdot(2AB^{2}+B^{3}d)+0\\
+((A^{4}d+2A^{3}B)d+0)\cdot(A^{3}d+3A^{2}B+AB^{2}d)\\
+(2A^{3}Bd+A^{2}B^{2}(5+d^{2})+4AB^{3}d+B^{4}d^{2})\cdot((2AB^{2}%
+B^{3}d)d+0)\\
+(2A^{3}Bd+A^{2}B^{2}(5+d^{2})+4AB^{3}d+B^{4}d^{2})\cdot(A^{3}d+3A^{2}%
B+AB^{2}d).
\end{gather*}
The writhe of $D$ is 3, so we calculate $f_{D}(A)$ by multiplying by $-A^{-9}
$ and evaluating $d\mapsto-A^{2}-B^{2}$ and $B\mapsto A^{-1}$:
\begin{align*}
f_{D}(A)  & =-A^{-9}\cdot(-A^{6}+A^{2})\cdot(A^{-1}-A^{-5})\\
& -A^{-9}\cdot(-A^{6}+A^{2})(-A^{2}-A^{-2})\cdot(-A^{5}+A-A^{-3})\\
& -A^{-9}\cdot(-A^{4}+2-A^{-4}+A^{-8})\cdot(-A^{-5}+A^{-1})(-A^{2}-A^{-2})\\
& -A^{-9}\cdot(-A^{4}+2-A^{-4}+A^{-8})\cdot(-A^{5}+A-A^{-3})\\
& =-A^{-9}\cdot(-A^{13}+2A^{9}-3A^{5}+3A-4A^{-3}+3A^{-7}-2A^{-11}+A^{-15})
\end{align*}
We conclude that the Jones polynomial of the knot of Figure \ref{tanglex} is
$f_{D}(t^{-1/4})=-t^{6}+2t^{5}-3t^{4}+4t^{3}-3t^{2}+3t-2+t^{-1}$. This
identifies the knot as the mirror image of $7_{6}$.

\end{document}